# HOEFFDING'S LEMMA FOR GENERAL MARKOV CHAINS AND ITS APPLICATIONS TO STATISTICAL LEARNING

By Jianqing Fan[‡], Bai Jiang[‡,†] and Qiang Sun[§]

*Princeton University*[‡], *University of Toronto*[§]

We extend Hoeffding's lemma to general-state-space and not necessarily reversible Markov chains. Let $\{X_i\}_{i\geq 1}$ be a stationary Markov chain with invariant measure $\pi$ and absolute spectral gap $1 - \lambda \in [0, 1]$, where $\lambda$ is defined as the operator norm of the transition kernel acting on mean zero and square-integrable functions with respect to $\pi$. Then, for any bounded functions $f_i : x \mapsto [a_i, b_i]$ and any $t \in \mathbb{R}$,

$$\mathbb{E}\left[e^{t \sum_i [f_i(X_i) - \int f_i(x)\pi(dx)]}\right] \leq \exp\left(\frac{t^2}{2} \cdot \frac{1+\lambda}{1-\lambda} \cdot \sum_i \frac{(b_i - a_i)^2}{4}\right).$$

This bound differs from the classical Hoeffding's lemma by a multiplicative coefficient of $(1 + \lambda)/(1 - \lambda)$, and simplifies to the latter when $\lambda = 0$. The counterpart of Hoeffding's inequality for Markov chains immediately follows. Our results assume none of countable state space, reversibility and time-homogeneity of Markov chains and cover time-dependent functions with various ranges. We illustrate the utility of these results by applying them to six problems in statistics and machine learning.

## 1. Introduction.

1.1. *Hoeffding's inequality for Markov chains.* Hoeffding (1963) proved in his celebrated work that a bounded random variable $Z \in [a, b]$ is sub-Gaussian[1] with variance proxy $(b-a)^2/4$. It follows that the sum of $n$ independent, bounded random variables $Z_i \in [a_i, b_i]$, $i = 1, \ldots, n$ is sub-Gaussian with variance proxy $\sum_{i=1}^{n}(b_i - a_i)^2/4$. Specifically, for any $t \in \mathbb{R}$,

$$\mathbb{E}\left[e^{t \sum_{i=1}^n (Z_i - \mathbb{E} Z_i)}\right] \leq \exp\left(\frac{t^2}{2} \cdot \sum_{i=1}^n \frac{(b_i - a_i)^2}{4}\right). \quad (1.1)$$

From this bound, the Chernoff approach (Boucheron, Lugosi and Massart, 2013) derives Hoeffding's inequality, which controls the tail probability for the sum of

---

*Supported by NSF under grant numbers DMS-1662139 and DMS-1712591 and NIH under grant number 2R01-GM072611-12. The bulk of this work was completed while Qiang Sun was a postdoc researcher at Princeton University.

[†]To whom correspondence should be addressed.

*MSC 2010 subject classifications:* Primary 60J05; secondary 62J05, 65C05

*Keywords and phrases:* general Markov chain, Hoeffding's lemma, Hoeffding's inequality, statistical learning

[1]A random variable $Z$ is sub-Gaussian with variance proxy $\sigma^2$ if $Z$ has a finite mean $\mathbb{E}Z$ and $\mathbb{E}[e^{t(Z-\mathbb{E}Z)}] \leq \exp(\sigma^2 t^2/2)$ for any $t \in \mathbb{R}$. An $n$-dimensional random vector $\boldsymbol{Z}$ is sub-Gaussian with variance proxy $\sigma^2$ if $\boldsymbol{Z}$ has a finite mean $\mathbb{E}\boldsymbol{Z}$ and $\mathbb{E}[e^{\boldsymbol{t}'(\boldsymbol{Z}-\mathbb{E}\boldsymbol{Z})}] \leq \exp(\sigma^2 \|\boldsymbol{t}\|_2^2/2)$ for any $\boldsymbol{t} \in \mathbb{R}^n$.

independent, bounded random variables. For any $\epsilon > 0$,

$$\mathbb{P}\left(\left|\sum_{i=1}^{n} Z_i - \sum_{i=1}^{n} \mathbb{E} Z_i\right| > \epsilon\right) \leq 2 \exp\left(-\frac{\epsilon^2}{2\sum_{i=1}^{n}(b_i - a_i)^2/4}\right). \qquad (1.2)$$

Various authors have discovered Hoeffding-type inequalities for Markov dependent random variables $Z_i = f(X_i)$, where $\{X_i\}_{i \geq 1}$ is a Markov chain and $f : x \mapsto [a, b]$ is a bounded function, by spectral methods. Their inequalities involve spectral gaps as coefficients. Among them, Gillman (1993) obtained the first Hoeffding-type inequality for finite-state-space and reversible Markov chains, which was later improved by Dinwoodie (1995). In the same setting, León and Perron (2004) provided a sharp Hoeffding-type inequality which differed from the classical Hoeffding's inequality by a multiplicative coefficient. Miasojedow (2014) extended León and Perron's Hoeffding-type inequality to general-state-space and not necessarily reversible Markov chains, but his inequality has a looser multiplicative coefficient than León and Perron's. Chung et al. (2012) independently established another interesting Hoeffding-type inequality for finite-state-space and not necessarily reversible Markov chains.

1.2. *Summary of our results.* Following (León and Perron, 2004) and (Miasojedow, 2014), this paper establishes the exact counterparts of the classical Hoeffding's lemma (1.1) and inequality (1.2) for general-state-space and not necessarily reversible Markov chains. For convenience, we refer to these Markov chains as general Markov chains.

Let us introduce some notations before presenting our main results. Denote by $\pi$ the invariant measure of the Markov chain $\{X_i\}_{i \geq 1}$. Let $\mathcal{L}_2(\pi)$ be the Hilbert space consisting of square-integrable functions with respect to $\pi$, and $\mathcal{L}_2^0(\pi)$ be its subspace of mean zero functions. Denote $P$ as the transition kernel of the Markov chain, which is viewed as a Markov operator acting on $\mathcal{L}_2(\pi)$, and by $P^*$ its adjoint. Let $\lambda \in [0, 1]$ be the operator norm of $P$ acting on $\mathcal{L}_2^0(\pi)$. We refer to $1 - \lambda$ as the *absolute spectral gap* of the Markov chain. This quantity measures the converging speed of the Markov chain toward its invariant measure $\pi$ (Rudolf, 2012).

Following Fill (1991), we refer to $(P + P^*)/2$ as the additive reversiblization of $P$, which is also a Markov operator acting on $\mathcal{L}_2(\pi)$. For a reversible Markov chain, $P$ is self-adjoint and coincides with $P^*$ and $(P + P^*)/2$. For a non-reversible Markov chain, $P$ is not self-adjoint, but $(P + P^*)/2$ is self-adjoint and associates with a reversible transition kernel. It is known that the spectrum of a self-adjoint Markov operator, such as $(P + P^*)/2$, acting on $\mathcal{L}_2^0(\pi)$ is contained in $[-1, +1]$ on the real line. We let $\lambda_r \in [-1, +1]$ denote the rightmost value of the spectrum of $(P + P^*)/2$ acting on $\mathcal{L}_2^0(\pi)$, and refer to $1 - \lambda_r$ as the *right spectral gap* of the Markov chain. Note that $\lambda_r$ is smaller than or equal to the operator norm of $(P + P^*)/2$ acting on $\mathcal{L}_2^0(\pi)$, which is then smaller than or equal to that of $P$, and thus $\lambda_r \leq \lambda$.

Our first main result, stated as Theorem 2.1, considers a stationary Markov chain $\{X_i\}_{i \geq 1}$ and time-dependent bounded functions $f_i : x \mapsto [a_i, b_i]$ with various ranges, and asserts that $\sum_{i=1}^{n} f_i(X_i)$ is sub-Gaussian with variance proxy

$$\alpha(\lambda) \times \sum_{i=1}^{n}(b_i - a_i)^2/4, \text{ where } \alpha : \lambda \mapsto \frac{1+\lambda}{1-\lambda} \geq 1.$$

This theorem and its resulting Hoeffding's inequality simplify to the classical Hoeffding's lemma (1.1) and inequality (1.2) when $\lambda = 0$. The resulting Hoeffding's inequality generalizes Miasojedow (2014)'s inequality from the time-independent function case, in which $f_1 = f_2 = \cdots = f_n = f$ are identical, to the time-dependent function case, in which $f_i$ are not identical.

For the time-independent function case in which $f_1 = f_2 = \cdots = f_n = f$ are identical, we find that it is possible to sharpen the multiplicative coefficient in Miasojedow (2014)'s and our Theorem 2.1 by replacing $\lambda$ with a smaller quantity $\max\{\lambda_r, 0\}$. This result, formally stated in Theorem 2.2, precisely generalizes León and Perron (2004)'s inequality for finite-state-space and reversible Markov chains to that for general Markov chains. Note that, for such Markov chains, $\lambda_r$ is merely the second largest eigenvalue of the transition probability matrix.

It is also worth noting that the our theorems discover the sub-Gaussian property of $\sum_{i=1}^{n} f_i(X_i)$ or $\sum_{i=1}^{n} f(X_i)$, which had been missing in the previous literature. This sub-Gaussian property, we believe, is the essence of Hoeffding's inequality and explicitly characterizes the behaviors of Markov dependent random variables.

1.3. *Other related literature.* In addition to Hoeffding-type inequalities, spectral methods have produced Bernstein-type inequalities, see e.g. (Lezaud, 1998a, Theorem 1.1) and (Paulin, 2015, Theorems 3.3 and 3.4). For finite-state-space Markov chains, their method consists of three steps: first bound the moment generating function (mgf) of $\sum_i f_i(X_i)$ or $\sum_i f(X_i)$ by a product of norms of perturbed Markov operators, then use the fact that the norms of operators on finite state spaces coincide with their largest eigenvalues, finally estimate their largest eigenvalues by Kato (2013)'s perturbation theory. Lezaud (1998b, 2001) attempted to extend this method to general state spaces and showed that the final step holds in the general setting. However, the second step does not follow easily because an operator in the general setting is infinite-dimensional, has complicated spectrum, and thus its norm may not coincide with its largest eigenvalue. After attentively investigating into the related works, we confirm that the proofs of (Lezaud, 1998a, Theorem 1.1) and (Paulin, 2015, Theorems 3.3 and 3.4) in general state spaces are incomplete.

Other exponential concentration inequalities for the sum of $n$ functions exploit the minorization and drift conditions, see e.g. (Glynn and Ormoneit, 2002; Douc et al., 2011; Adamczak and Bednorz, 2015), or information-theoretical ideas, see e.g. (Kontoyiannis et al., 2005, 2006). Their inequalities apply for Markov chains which may not admit non-zero spectral gaps, but have less explicit and often suboptimal constants or more complicated expressions.

There is another and less related research line for a function of $n$ variables under Markov or other dependent structures. Marton (1996, 1998) pioneered the concentration of measure phenomenon for contracting Markov chains. Further progresses are made by Samson (2000); Chazottes et al. (2007); Redig and Chazottes (2009); Kontorovich and Ramanan (2008), among many others.

1.4. *Applications to statistical learning problems.* As well as the theoretical interest, this work is motivated by the practical needs of concentration inequalities for Markov chains arising in modern statistics and machine learning problems.

First, we consider linear regression, lasso regression and sparse covariance matrix estimation for time series data with Markovian structures, see Yuan and Kendziorski (2006); Nardi and Rinaldo (2011); Cai et al. (2010); Fan et al. (2016) for examples. Existing analyses of these methods usually assume independent data samples. In order to close the gap between theory and practice, we apply the derived concentration inequalities for Markov-dependent random variables to obtain optimal non-asymptotic error bounds for estimators.

Another motivating example is the Markov chain Monte Carlo (MCMC) (Gilks, 2005) approximation method, which evaluates a complicated integral by averaging samples from a well-behaved Markov chain. A sharp concentration inequality is needed to determine how long the Markov chain should run in order to control the approximation error.

Our last example is the multi-armed bandit problem with Markovian rewards. Researchers have recognized three fundamental formalizations of multi-armed bandit problems depending on the nature of the reward process: independently and identically distributed (i.i.d.), adversarial, and Markovian (Bubeck and Cesa-Bianchi, 2012). The Markovian formalization has been much less studied in the way of the other two, primarily because bandit algorithms rely on the concentration of reward draws to identify the optimal arm, but unfortunately powerful concentration inequalities for Markovian rewards are sparse.

1.5. *Organization of Paper.* The rest of this paper proceeds as follows. Section 2 presents our central results for a stationary, homogeneous Markov chain in Theorem 2.1 (for time-dependent functions) and Theorem 2.2 (for a time-independent function), and discusses its extensions to non-stationary, non-homogeneous Markov chains and unbounded functions. An impossibility result is given for unbounded functions. Section 3 introduces preliminaries of Hilbert spaces and the operator theory. Section 4 structures the proof of Theorem 2.1. Section 5 is devoted to the proof of Theorem 2.2. Section 6 applies our theorems to six problems in statistics and machine learning.

**2. Main Results.** Let $\mathcal{X}$ denote the general state space of the Markov chain $\{X_i\}_{i\geq 1}$ under study. Let $\mathbb{E}_\nu[\cdot]$ and $\mathbb{P}_\nu(\cdot)$ denote the expectation of random variables and the probability of events when the Markov chain starts from an initial measure $\nu$. Recall that $\alpha : \lambda \mapsto (1+\lambda)/(1-\lambda)$. For any function $h : \mathcal{X} \to \mathbb{R}$, we write $\pi(h) \coloneqq \int h(x)\pi(dx)$ whenever the integral is finite. Write $\xi_1 \vee \xi_2 = \max\{\xi_1, \xi_2\}$ for $\xi_1, \xi_2 \in \mathbb{R}$.

Our first central result is formally stated as the following theorem.

THEOREM 2.1. *Let $\{X_i\}_{i\geq 1}$ be a Markov chain with invariant measure $\pi$ and absolute spectral gap $1-\lambda > 0$ (see Definition 3.6). Then for any bounded functions $f_i : \mathcal{X} \to [a_i, b_i]$, the sum $\sum_{i=1}^n f_i(X_i)$ is sub-Gaussian with variance proxy $\alpha(\lambda) \cdot \sum_{i=1}^n (b_i - a_i)^2/4$. That is, for any $t \in \mathbb{R}$*

$$\mathbb{E}_\pi \left[ e^{t \sum_{i=1}^n (f_i(X_i) - \pi(f_i))} \right] \leq \exp\left( \frac{t^2}{2} \cdot \alpha(\lambda) \cdot \sum_{i=1}^n \frac{(b_i - a_i)^2}{4} \right). \quad (2.1)$$

*It follows that for any $\epsilon > 0$*

$$\mathbb{P}_\pi \left( \left| \sum_{i=1}^n f_i(X_i) - \sum_{i=1}^n \pi(f_i) \right| > \epsilon \right) \leq 2 \exp\left( -\frac{\alpha(\lambda)^{-1} \epsilon^2}{2 \sum_{i=1}^n (b_i - a_i)^2/4} \right). \quad (2.2)$$

As stated in (Rudolf, 2012), $1 - \lambda$ quantifies the converging speed of the Markov chain toward its invariant distribution $\pi$. A smaller $\lambda$ indicates faster convergence speed and less variable dependence. The multiplicative coefficient $\alpha(\lambda) = (1 + \lambda)/(1 - \lambda)$ in Theorem 2.1 is strictly increasing with $\lambda$. This is consistent with the intuition that a faster Markov chain $\{X_i\}_{i\geq 1}$ with less dependence among variables results in a smaller variance proxy of $\sum_{i=1}^n f_i(X_i)$, and thus a tighter concentration.

When $\lambda = 0$ and correspondingly $\alpha(\lambda) = 1$, (2.1) and (2.2) in Theorem 2.1 reduce to the classical Hoeffding's lemma (1.1) and inequality (1.2) for independent random variables. Indeed, independent random variables $Z_i \in [a_i, b_i]$ can be seen as transformations of independently and identically distributed (i.i.d.) random variables $U_i \sim \text{Uniform}[0, 1]$ via the inverse cumulative distribution functions $F_{Z_i}^{-1} : [0, 1] \to [a_i, b_i]$. Namely, $Z_i = F_{Z_i}^{-1}(U_i)$. The sequence of i.i.d. $\{U_i\}_{i\geq 1}$ forms a stationary Markov chain on the state space $[0, 1]$ with invariant measure $\pi(dy) = dy$ and transition kernel $P(x, dy) = dy$. This Markov chain has $\lambda = 0$.

We also find that $\alpha(\lambda)$ is the smallest multiplicative coefficient for Theorem 2.1 to hold. Consider the case in which $f_1 = \ldots f_n = f$ are Rademacher, i.e. $\pi(\{f = +1\}) = \pi(\{f = -1\}) = 1/2$, and the transition kernel is $P(x, B) = \lambda \mathbb{I}(x \in B) + (1-\lambda)\pi(B)$ for any state $x \in \mathcal{X}$ and any subset of the state space $B \subseteq \mathcal{X}$. Applying (Geyer, 1992, Theorem 2.1) yields

$$\lim_{n\to\infty} \text{Var}\left( \frac{1}{\sqrt{n}} \sum_{i=1}^n f(X_i) \right) = \alpha(\lambda).$$

As the variance proxy of a sub-Gaussian random variable naturally upper bounds its variance, the variance proxy of $\sum_{i=1}^n f(X_i)$ cannot be smaller than $n\alpha(\lambda)$. This lower bound is achieved by Theorem 2.1.

2.1. *Result for time-independent functions.* For the time-independent function case in which $f_1 = \cdots = f_n = f$ are identical, Theorem 2.1 coincides with the Hoeffding-type inequality in (Miasojedow, 2014). Further, we find that it is possible to sharpen the multiplicative coefficient $\alpha(\lambda)$ in Miasojedow's inequality and our Theorem 2.1 by substituting the absolute spectral gap $1 - \lambda$ with the right spectral gap $1 - \lambda_r$. This sharper inequality is presented in Theorem 2.2. Recall that $\xi_1 \vee \xi_2 = \max\{\xi_1, \xi_2\}$ for $\xi_1, \xi_2 \in \mathbb{R}$.

THEOREM 2.2. *Let $\{X_i\}_{i\geq 1}$ be a Markov chain with invariant measure $\pi$ and right spectral gap $1 - \lambda_r > 0$ (see Definition 3.7). Then for any bounded function $f : \mathcal{X} \to [a, b]$ and any $t \in \mathbb{R}$,*

$$\mathbb{E}_\pi \left[ e^{t\left(\sum_{i=1}^n f(X_i) - n\pi(f)\right)} \right] \leq \exp\left( \frac{t^2}{2} \cdot \alpha(\lambda_r \vee 0) \cdot \frac{n(b-a)^2}{4} \right). \qquad (2.3)$$

*It follows that for any $\epsilon > 0$,*

$$\mathbb{P}_\pi \left( \left| \sum_{i=1}^n f(X_i) - n\pi(f) \right| > \epsilon \right) \leq 2 \exp\left( -\frac{\alpha(\lambda_r \vee 0)^{-1} \epsilon^2}{2n(b-a)^2/4} \right). \qquad (2.4)$$

Recall that $\lambda \geq \lambda_r$ and $\lambda \geq 0$ (see the preliminary section for more technical details) and that $\alpha : \lambda \mapsto (1 + \lambda)/(1 - \lambda)$ is a strictly increasing function. The multiplicative coefficient $\alpha(\lambda_r \vee 0)$ in the proceeding theorem is sharper than the previous multiplicative coefficient $\alpha(\lambda)$ in (Miasojedow, 2014) and Theorem 2.1. Also, the condition $\lambda_r < 1$ in the former is weaker than the condition $\lambda < 1$ in the latter. Moreover, for finite-state-space and reversible Markov chains studied by León and Perron (2004), $\lambda_r$ coincides with the second largest eigenvalue of the transition probability matrix. For these Markov chains, the Hoeffding-type inequality (2.4) precisely reduces to León and Perron (2004)'s, whereas Miasojedow's and Theorem 2.1 do not.

2.2. *Result for non-stationary Markov chains.* Theorem 2.1 and Theorem 2.2 are concerned with stationary Markov chains. For a non-stationary Markov chain starting from an initial measure $\nu \neq \pi$, Theorem 2.3 presents similar inequalities given that the initial measure $\nu$ satisfies some moment condition. The proof is put in Subsection A.1 in the appendix.

THEOREM 2.3. *Let $\{X_i\}_{i\geq 1}$ be a Markov chain with invariant measure $\pi$ and absolute spectral gap $1 - \lambda > 0$ (see Definition 3.6). Suppose the initial measure $\nu$ is absolutely continuous with respect to the invariant measure $\pi$ and its density, denoted by $d\nu/d\pi$, has a finite p-moment for some $p \in (1, \infty]$, i.e.*

$$\infty > \left\| \frac{d\nu}{d\pi} \right\|_{\pi,p} := \begin{cases} \left[ \pi\left( \left|\frac{d\nu}{d\pi}\right|^p \right) \right]^{1/p} & \text{if } p < \infty, \\ \text{ess sup} \left| \frac{d\nu}{d\pi} \right| & \text{if } p = \infty. \end{cases}$$

*Let $q = p/(p-1) \in [1, \infty)$. Then for any bounded functions $f_i : \mathcal{X} \to [a_i, b_i]$ and any $t \in \mathbb{R}$,*

$$\mathbb{E}_\nu \left[ e^{t \sum_{i=1}^n (f_i(X_i) - \pi(f_i))} \right] \leq \left\| \frac{d\nu}{d\pi} \right\|_{\pi,p} \exp\left( \frac{t^2}{2} \cdot q\alpha(\lambda) \cdot \sum_{i=1}^n \frac{(b_i - a_i)^2}{4} \right).$$

*It follows that for any $\epsilon > 0$*

$$\mathbb{P}_\nu \left( \left| \sum_{i=1}^n f_i(X_i) - \sum_{i=1}^n \pi(f_i) \right| > \epsilon \right) \leq 2 \left\| \frac{d\nu}{d\pi} \right\|_{\pi,p} \exp\left( -\frac{q^{-1}\alpha(\lambda)^{-1} \epsilon^2}{2\sum_{i=1}^n (b_i - a_i)^2/4} \right).$$

The sub-Gaussian variance proxy for the non-stationary Markov chains is $q = p/(p-1) \geq 1$ times that for the stationary Markov chains. This larger variance proxy leads to a slower concentration in the Hoeffding-type inequality. If $d\nu/d\pi$ is essentially bounded by some absolute constant then $p = \infty$, $q = 1$, $\|d\nu/d\pi\|_{\pi,\infty} < \infty$ obtaining the same exponent rate in the Hoeffding-type bound with the stationary case. We later show in Subsection 6.4 that a combination of Theorems 2.2 and 2.3 provides a non-asymptotic error bound for MCMC integrals.

2.3. *Results for unbounded (sub-Gaussian) functions.* Surprisingly, the assumption on the boundedness of functions in Theorems 2.1 and 2.2 cannot be relaxed to sub-Gaussian functions. Below is an example of a stationary Markov chain and a time-independent function $f_i = f$, in which $f(X_i)$ individually follows the standard Gaussian distribution $\mathcal{N}(0,1)$ but no finite multiplicative coefficient exists for (2.1), (2.2), (2.3) and (2.4) to hold. The proof is collected in Subection A.2 in the appendix.

THEOREM 2.4. *Consider a stationary Markov chain $\{X_i\}_{i\geq 0}$ on state space $\mathcal{X} = \mathbb{R}$ with invariant distribution $\pi \sim \mathcal{N}(0,1)$ and transition kernel*

$$P(x, B) = \lambda \mathbb{I}(x \in B) + (1-\lambda)\pi(B), \quad \forall x \in \mathcal{X}, \ \forall B \subseteq \mathcal{X}.$$

*There exists no finite constant $\widetilde{\alpha}$ such that $\sum_{i=1}^n X_i$, $n \geq 1$ are sub-Gaussian with variance proxy $n\widetilde{\alpha}$.*

2.4. *Results for time-inhomogeneous Markov chains.* We also consider relaxing the implicit assumption on the time-homogeneity of Markov chains in Theorem 2.1. For a time-inhomogeneous Markov chain with different transition kernels $\{P_i\}_{i\geq 1}$ at each step $i$, we write $1 - \lambda_i$ as the absolute spectral gap of $P_i$. A refinement of the proof of Theorem 2.1 shows that $\sum_{i=1}^n f_i(X_i)$ is sub-Gaussian with variance proxy

$$\frac{(b_1 - a_1)^2}{8} + \sum_{i=2}^n \alpha(\lambda_{i-1}) \cdot \frac{(b_{i-1} - a_{i-1})^2 + (b_i - a_i)^2}{8} + \frac{(b_n - a_n)^2}{8}$$
$$\leq \alpha\left(\max_{1\leq i\leq n} \lambda_i\right) \cdot \sum_{i=1}^n \frac{(b_i - a_i)^2}{4}.$$

If time-inhomogeneous transition kernels admit a uniform absolute spectral gap $1 - \lambda$, i.e. $\lambda_i \leq \lambda$ for all $i \geq 1$, then we have (2.1) and (2.2). This result cannot be established by arguments in León and Perron (2004) and Miasojedow (2014).

**3. Preliminaries.** Throughout the paper, we assume the state space $\mathcal{X}$ equipped with a sigma-algebra $\mathcal{B}$ is *standard Borel*[2]. This is a common assumption to rigorously define and investigate Markov chains in measure theory. In most practical

---

[2] A measurable space $(\mathcal{X}, \mathcal{B})$ is standard Borel if it is isomorphic to a subset of $\mathbb{R}$. Such measurable spaces are also called *nice* spaces. See (Breiman, 1992, Definition 4.33) and (Durrett, 2010, page 45)

examples, $\mathcal{X}$ is a subset of a multi-dimensional real space and $\mathcal{B}$ is the Borel sigma-algebra over $\mathcal{X}$.

The distribution of a time-homogeneous Markov chain is uniquely determined by its initial measure $\nu$ and its transition kernel $P$.

$$\nu(B) = \mathbb{P}(X_1 \in B), \quad \forall B \in \mathcal{B}.$$
$$P(X_i, B) = \mathbb{P}(X_{i+1} \in B | X_i), \quad \forall B \in \mathcal{B}, \ \forall i \geq 0.$$

A transition kernel $P$ is invariant with a probability measure $\pi$ on $(\mathcal{X}, \mathcal{B})$ if

$$\pi(B) = \int P(x, B)\pi(dx), \quad \forall B \in \mathcal{B}.$$

A Markov chain is said stationary if it starts from $\nu = \pi$.

Our analyses are conducted in the framework of operator theory on Hilbert spaces. The idea of using this framework originates from the fact that each transition kernel, if invariant with $\pi$, is viewed as a Markov operator on the Hilbert space $\mathcal{L}_2(\pi)$ consisting of all real-valued, $\mathcal{B}$-measurable, $\pi$-square-integrable functions on $\mathcal{X}$.

3.1. *Hilbert Space $\mathcal{L}_2(\pi)$*. Recall that we write $\pi(h) := \int h(x)\pi(dx)$ for any real-valued, $\mathcal{B}$-measurable function $h : \mathcal{X} \to \mathbb{R}$. Let $\mathcal{L}_p(\mathcal{X}, \mathcal{B}, \pi)$ be the set of real-valued, $\mathcal{B}$-measurable functions with finite $p$-moment, i.e.

$$\mathcal{L}_p(\mathcal{X}, \mathcal{B}, \pi) := \{h : \pi(|h|^p) < \infty\}.$$

Here $h_1, h_2 \in \mathcal{L}_p(\mathcal{X}, \mathcal{B}, \pi)$ are taken as identical if $h_1 = h_2$ $\pi$-almost everywhere ($\pi$-a.e.). For every $p \in [1, \infty]$, $\mathcal{L}_p(\mathcal{X}, \mathcal{B}, \pi)$ is a Banach space equipped with norm

$$\|h\|_{\pi,p} := \begin{cases} \pi(|h|^p)^{1/p} & \text{if } p < \infty, \\ \operatorname{ess\,sup} |h| & \text{if } p = \infty. \end{cases}$$

In particular, if $p = 2$ then

$$\mathcal{L}_2(\mathcal{X}, \mathcal{B}, \pi) := \{h : \pi(h^2) < \infty\}$$

is a Hilbert space[3] endowed with the following inner product

$$\langle h_1, h_2 \rangle_\pi = \int h_1(x) h_2(x) \pi(dx), \quad \forall h_1, h_2 \in \mathcal{L}_2(\mathcal{X}, \mathcal{B}, \pi),$$

since for every $h \in \mathcal{L}_2(\mathcal{X}, \mathcal{B}, \pi)$,

$$\|h\|_{\pi,2} = \sqrt{\langle h, h \rangle_\pi}.$$

By convention, the norm of a linear operator $T$ on $\mathcal{L}_2(\mathcal{X}, \mathcal{B}, \pi)$ is defined as

$$\|\!|T|\!\|_{\pi,2} = \sup\{\|Th\|_{\pi,2} : \|h\|_{\pi,2} = 1\}.$$

---

[3]Here we consider this real Hilbert space instead of the complex Hilbert space, as the former is adequate for our proofs.

Another important Hilbert space is the subspace of $\mathcal{L}_2(\mathcal{X}, \mathcal{B}, \pi)$ consisting of mean zero functions.

$$\mathcal{L}_2^0(\mathcal{X}, \mathcal{B}, \pi) := \{h \in \mathcal{L}_2(\mathcal{X}, \mathcal{B}, \pi) : \pi(h) = 0\}.$$

For simplicity of notations, we write $\|h\|_\pi$ and $\|\|T\|\|_\pi$ in place of $\|h\|_{\pi,2}$ and $\|\|T\|\|_{\pi,2}$, respectively. We also write $\mathcal{L}_2(\pi)$ and $\mathcal{L}_2^0(\pi)$ in place of $\mathcal{L}_2(\mathcal{X}, \mathcal{B}, \pi)$ and $\mathcal{L}_2^0(\mathcal{X}, \mathcal{B}, \pi)$ respectively, whenever the measurable space $(\mathcal{X}, \mathcal{B})$ is clear in the context.

3.2. *Viewing transition kernels as operators on $\mathcal{L}_2(\pi)$.* Each transition kernel $P(x, B)$, if invariant with $\pi$, corresponds to a bounded linear operator on $\mathcal{L}_2(\pi)$. We call this operator Markov and abuse $P$ to denote it. That is,

$$Ph(x) = \int h(y) P(x, dy), \quad \forall x \in \mathcal{X}, \ \forall h \in \mathcal{L}_2(\pi).$$

Next, we introduce five transition kernels and their associated Markov operators which appear frequently throughout the rest of the paper. We abuse the same notation for a transition kernel and its associated Markov operator.

3.2.1. *Identity operator $I$.* The identity kernel given by

$$I(x, B) = \mathbb{I}(x \in B), \quad \forall x \in \mathcal{X}, \ \forall B \in \mathcal{B},$$

generates a Markov chain which never moves from its initial state. The identity kernel corresponds to the identity operator on $\mathcal{L}_2(\pi)$

$$I : h \in \mathcal{L}_2(\pi) \mapsto h.$$

3.2.2. *Projection operator $\Pi$.* The transition kernel given by

$$\Pi(x, B) = \pi(B), \quad \forall x \in \mathcal{X}, \ \forall B \in \mathcal{B},$$

generates a Markov chain which consists of i.i.d. draws from the invariant measure $\pi$. Denote by the italicized symbol *1* the constant function $x \in \mathcal{X} \mapsto 1$. The transition kernel $\Pi(x, B)$ corresponds to the following Markov operator

$$\Pi : h \in \mathcal{L}_2(\pi) \mapsto \pi(h)\mathit{1},$$

which is a projection operator of rank one since $\pi(h) = \langle h, \mathit{1}\rangle_\pi$.

If the transition kernel $P(x, B)$ is invariant with $\pi$, then $P\Pi = \Pi P = \Pi$.

3.2.3. *Adjoint operator $P^*$.*

DEFINITION 3.1 (Time-reversal kernel). *A transition kernel $P^*$ is said to be the time-reversal of a transition kernel $P$ if*

$$\int_{B_1} \pi(dx) P(x, B_2) = \int_{B_2} \pi(dx) P^*(x, B_1), \quad \forall B_1, B_2 \in \mathcal{B}.$$

The existence of this time-reversal kernel in a standard Borel space $(\mathcal{X}, \mathcal{B})$ is assured by (Breiman, 1992, Theorem 4.34) or (Durrett, 2010, Theorem 5.1.9). It is also unique up to differences on sets of probability zero. In most practical examples, probability measures $\pi(\cdot)$ and $\{P(x, \cdot) : x \in \mathcal{X}\}$ share a common reference measure. In these examples, let $\pi(x)$ and $P(x,y)$ denote their densities then the time-reversal transition kernel has a density of a simple closed form

$$P^*(x,y) = \frac{\pi(y)P(y,x)}{\pi(x)}.$$

The time-reversal kernel corresponds to the adjoint operator (see Definition 3.2) of $P$ on $\mathcal{L}_2(\pi)$.

DEFINITION 3.2 (Adjoint operator). *A linear operator $T^*$ on a real Hilbert space $\mathcal{H}$ endowed with inner product $\langle \cdot, \cdot \rangle$ is said to be the adjoint of a linear operator $T$ if*

$$\langle Th_1, h_2 \rangle = \langle h_1, T^*h_2 \rangle, \quad \forall h_1, h_2 \in \mathcal{H}.$$

If the transition kernel $P(x, B)$ is invariant with $\pi$, so is $P^*(x, B)$. And, $P^*\Pi = \Pi P^* = \Pi$.

A Markov chain is reversible if $P^* = P$. This condition is called the detailed balance when viewing $P$ and $P^*$ as transition kernels, or the self-adjointness (see Definition 3.3) when viewing $P$ and $P^*$ as Markov operators.

DEFINITION 3.3 (Self-adjoint operator). *A linear operator $S$ on a real Hilbert space $\mathcal{H}$ endowed with inner product $\langle \cdot, \cdot \rangle$ is said to be self-adjoint if*

$$\langle Sh_1, h_2 \rangle = \langle h_1, Sh_2 \rangle, \quad \forall h_1, h_2 \in \mathcal{H}.$$

3.2.4. *Additively-reversiblized operator $R$.* Fill (1991) defined the additive reversiblization[4] of a Markov operator $P$ as

$$R = \frac{P + P^*}{2},$$

which is self-adjoint and thus relates to a reversible Markov transition kernel $R(x, B)$. If the transition kernel $P(x, B)$ is invariant with $\pi$, so is $R(x, B)$. And, $R\Pi = \Pi R = \Pi$.

3.2.5. *León-Perron Operators $\widehat{P}$.* Every convex combination of Markov transition kernels (operators) produces a Markov transition kernel (operator). We say a Markov operator *León-Perron* if it is a convex combination of the identity operator $I$ and the projection operator $\Pi$.

---

[4]Authoritative books of operator theories refer to $R = (P + P^*)/2$ as the real part of $P$. Here we follow Fill (1991) as we discuss this operator in the context of Markov chains rather general operator theory

DEFINITION 3.4 (León-Perron operator). *A Markov operator $\widehat{P}_c$ on $\mathcal{L}_2(\pi)$ is said León-Perron if it is a convex combination of operators $I$ and $\Pi$ with some coefficient $c \in [0,1]$, that is*

$$\widehat{P}_c = cI + (1-c)\Pi.$$

The associated transition kernel

$$\begin{aligned}\widehat{P}_c(x, B) &= cI(x,B) + (1-c)\Pi(x,B) \\ &= c\mathbb{I}(x \in B) + (1-c)\pi(B), \quad \forall x \in \mathcal{X}, \ \forall B \in \mathcal{B},\end{aligned}$$

characterizes a random-scan mechanism: the Markov chain either stays at the current state (with probability $c$) or samples a new state from $\pi$ (with probability $1-c$) at each step.

If a León-Perron operator $\widehat{P}$ shares the same absolute spectral gap (see Definition 3.6) and invariant measure $\pi$ with a Markov operator $P$ then we call it the León-Perron version of $P$.

DEFINITION 3.5 (León-Perron version). *For a Markov operator $P$ with invariant measure $\pi$ and absolute spectral gap $1 - \lambda$, we say $\widehat{P}_\lambda = \lambda I + (1-\lambda)\Pi$ is the León-Perron version of $P$. We will omit the subscript $\lambda$ and write $\widehat{P}$ in place of $\widehat{P}_\lambda$ if $\lambda$ is clear in the context.*

The key to prove the Hoeffding-type inequalities for finite-state-space and reversible Markov chains in (León and Perron, 2004) is the observation that a Markov chain driven by this operator $\widehat{P}$ is the extremal case of all Markov chains with the same spectral gap. That is why we call this type of operators León-Perron.

3.3. *Absolute spectral gap.* Our theorems quantify the convergence speed of the Markov chain by the absolute spectral gap. This quantity has played a central role in the literature studying Markov chains by spectral methods, e.g. the central limit theorem for Metropolis-Hastings algorithm (Geyer, 1992), concentration inequalities of Markov-dependent random variables (Gillman, 1993; Dinwoodie, 1995; Lezaud, 1998a; León and Perron, 2004; Miasojedow, 2014; Paulin, 2015), the mean squared error of the Markov chain Monte Carlo estimators (Rudolf, 2012), and approximate transition kernels in Markov chain Monte Carlo algorithms (Negrea and Rosenthal, 2017).

Recall that $\mathcal{L}_2^0(\pi) \coloneqq \{h \in \mathcal{L}_2(\pi) : \pi(h) = 0\}$ denotes the subspace of $\mathcal{L}_2(\pi)$ consisting of mean zero functions.

DEFINITION 3.6 (absolute spectral gap). *A Markov operator $P$ admits an absolute spectral gap $1 - \lambda(P)$ if*

$$\lambda(P) := \sup\left\{\frac{\|Ph\|_\pi}{\|h\|_\pi} : h \in \mathcal{L}_2^0(\pi), h \neq 0\right\} < 1.$$

*Another equivalent characterizations of $\lambda(P)$ are given by*

$$\lambda(P) = \sup\{\|Ph\|_\pi : \|h\|_\pi = 1, \pi(h) = 0\} < 1,$$

*and*

$$\lambda(P) = \|\|P - \Pi\|\|_\pi < 1.$$

It is elementary that $\lambda(I) = 1$, $\lambda(\Pi) = 0$, $\lambda(P) = \lambda(P^*) \geq \lambda((P + P^*)/2)$, and $\lambda(\widehat{P}_c) = c$.

We briefly comment on the connection of this absolute spectral gap $1 - \lambda(P)$ to spectral radius and geometric ergodicity. Some literatures such as (Kontoyiannis and Meyn, 2012) define spectral gap as the gap between 1 and the spectral radius of $P$ acting on $\mathcal{L}_2^0(\pi)$. Denote by $\lambda_\infty(P)$ the spectral radius. It is known that

$$\lambda_\infty(P) = \lim_{k \to \infty} \|\|P^k - \Pi\|\|_\pi^{1/k} \leq \lambda(P).$$

The equality holds for reversible Markov chains. For these Markov chains, $\lambda_\infty(P) = \lambda(P) < 1$ if and only if the Markov chain is geometrically ergodic, see e.g. (Roberts and Rosenthal, 1997, Theorem 2.1) and (Conway, 2013, Proposition VIII.1.11). For non-reversible Markov chains, it is possible that $\lambda(P) > \lambda_\infty(P)$. Thereafter, the condition $1 - \lambda(P) > 0$ implies that $1 - \lambda_\infty(P) > 0$, and the latter further implies geometric ergodicity but not vice versa (Kontoyiannis and Meyn, 2012, Theorem 1.4). Under the weaker condition $1 - \lambda_\infty(P) > 0$, Theorem A.1 in the appendix give a similar Hoeffding-type inequality to Theorem 2.1.

3.4. *Right spectral gap.* Let $\mathcal{S}(R|\mathcal{L}_2^0(\pi))$ be the spectrum of the additive reversiblization

$$R = (P + P^*)/2$$

acting on $\mathcal{L}_2^0(\pi)$. It is known that the spectrum of such self-adjoint Markov operator on $\mathcal{L}_2^0(\pi)$ is contained in $[-1, +1]$ on the real line. Let

$$\begin{aligned}\lambda_{\rm r}(R) &:= \sup\{s : s \in \mathcal{S}(R|\mathcal{L}_2^0(\pi))\}, \\ \lambda_{\rm l}(R) &:= \inf\{s : s \in \mathcal{S}(R|\mathcal{L}_2^0(\pi))\}.\end{aligned} \quad (3.1)$$

Equivalent characterizations are given by

$$\begin{aligned}\lambda_{\rm r}(R) &:= \sup\{\langle Rh, h\rangle : \|h\|_\pi = 1, \ \pi(h) = 0\}, \\ \lambda_{\rm l}(R) &:= \inf\{\langle Rh, h\rangle : \|h\|_\pi = 1, \ \pi(h) = 0\}.\end{aligned}$$

Since $R$ is self-adjoint, it is known that

$$\lambda(R) = \sup\{|s| : s \in \mathcal{S}(R|\mathcal{L}_2^0(\pi))\} = \lambda_{\rm r}(R) \vee |\lambda_{\rm l}(R)|.$$

Thus

$$\lambda_{\rm r}(R) \leq \lambda(R) \leq \lambda(P).$$

DEFINITION 3.7 (Right spectral gap). *A Markov operator $P$ admits a right spectral gap $1 - \lambda_r$ if its additive reversiblization $R = (P + P^*)/2$ has $\lambda_{\rm r}(R) = \lambda_r < 1$.*

3.5. *Expressing the mgf of $\sum_i f_i(X_i)$ in $\mathcal{L}_2(\pi)$.* The essence of Theorem 2.1 is an upper bound for the mgf of $\sum_i f_i(X_i)$. In the Hilbert space $\mathcal{L}_2(\pi)$, this mgf has a simple expression involving the Markov operator $P$ and the multiplication operators of functions $e^{tf_i} : x \in \mathcal{X} \mapsto e^{tf_i(x)}$.

DEFINITION 3.8 (Multiplication Operator). *The multiplication operator of a function $g : \mathcal{X} \to \mathbb{R}$, denoted by $M_g$, is defined as*

$$M_g h(x) = g(x)h(x), \quad \forall x \in \mathcal{X}, \ \forall h \in \mathcal{L}_2(\pi).$$

$M_g$ is a bounded linear operator on $\mathcal{L}_2(\pi)$ if $g$ is bounded. For any bounded function $f$ and $t \in \mathbb{R}$, the function $e^{tf} : x \in \mathcal{X} \mapsto e^{tf(x)}$ is bounded and thus its multiplication operator $M_{e^{tf}}$ is bounded on $\mathcal{L}_2(\pi)$. To simplify the notations, we write $E^{tf}$ in place of $M_{e^{tf}}$.

Note that $E^{tf} = (E^{tf/2})^2$, that $E^{tf/2}$ is self-adjoint, and that $E^{tf} 1 = e^{tf}$. An elementary calculus expresses the mgf of $\sum_{i=1}^n f_i(X_i)$ as follows.

$$\begin{aligned}
\mathbb{E}_\pi \left[ e^{t \sum_{i=1}^n f_i(X_i)} \right] &= \left\langle 1, E^{tf_1} \left( \prod_{i=2}^n P E^{tf_i} \right) 1 \right\rangle_\pi \\
&= \left\langle 1, E^{tf_1/2} \left( \prod_{i=1}^{n-1} E^{tf_i/2} P E^{tf_{i+1}/2} \right) E^{tf_n/2} 1 \right\rangle_\pi \quad (3.2) \\
&= \left\langle e^{tf_1/2}, \left( \prod_{i=1}^{n-1} E^{tf_i/2} P E^{tf_{i+1}/2} \right) e^{tf_n/2} \right\rangle_\pi.
\end{aligned}$$

This expression will be the starting point of our proof.

**4. Proof of Theorem 2.1.** As defined above, $\widehat{P}_c$ denotes a general León-Perron operator with convex coefficient $c$, and the shorthand $\widehat{P} = \widehat{P}_{\lambda(P)}$ denotes the León-Perron version of $P$. Let $\{\widehat{X}_i\}_{i \geq 1}$ denote a Markov chain on the state space $\mathcal{X}$ driven by $\widehat{P}_c$ or $\widehat{P}$.

The proof of Theorem 2.1 proceeds as four subsections. Subsection 4.1 derives that for any $t \in \mathbb{R}$,

$$\mathbb{E}_\pi \left[ e^{t \sum_{i=1}^n f_i(X_i)} \right] \leq \prod_{i=1}^n \||E^{tf_i/2} \widehat{P} E^{tf_i/2}\||_\pi. \quad (4.1)$$

Thereafter it suffices to consider $\||E^{tf/2} \widehat{P} E^{tf/2}\||_\pi$ for function $f : \mathcal{X} \to [a,b]$ and the León-Perron version $\widehat{P}$ of $P$.

Subsection 4.2 shows that

$$\log \||E^{tf/2} \widehat{P} E^{tf/2}\||_\pi = \lim_{n \to \infty} \frac{1}{n} \log \mathbb{E}_\pi \left[ e^{t \sum_{i=1}^n f(\widehat{X}_i)} \right]. \quad (4.2)$$

Subsection 4.3 follows (4.2) to reduce the task of bounding $\||E^{tf/2} \widehat{P} E^{tf/2}\||_\pi$ into a simple problem in a two-state Markov chain system.

DEFINITION 4.1 (two-state Markov chain). *Let $\{\widehat{Y}_i\}_{i\geq 1}$ be a two-state Markov chain on the state space $\mathcal{Y} = \{a,b\}$ with a transition probability matrix $\widehat{\boldsymbol{Q}}$ determined by two parameters $\lambda \in [0,1)$ and $\mu \in (0,1)$ in the way*

$$\widehat{\boldsymbol{Q}} = \lambda \boldsymbol{I} + (1-\lambda) \begin{bmatrix} \boldsymbol{\mu}' \\ \boldsymbol{\mu}' \end{bmatrix}, \text{ where } \boldsymbol{\mu} = [1-\mu, \mu]'.$$

Note that $\widehat{\boldsymbol{Q}}$ is its own León-Perron version. Let $\boldsymbol{E}^{ty/2} = \mathrm{diag}(e^{ta/2}, e^{tb/2})$ be a $2\times 2$ diagonal matrix, and $\|\|\boldsymbol{T}\|\|_{\boldsymbol{\mu}}$ be the matrix norm of $\boldsymbol{T} \in \mathbb{R}^{2\times 2}$ induced by the $\boldsymbol{\mu}$-weighted vector norm. This subsection concludes that if $\mu = (\pi(f) - a)/(b-a)$ then for any $t \in \mathbb{R}$,

$$\|\|E^{tf/2}\widehat{P}E^{tf/2}\|\|_\pi \leq \|\|\boldsymbol{E}^{ty/2}\widehat{\boldsymbol{Q}}\boldsymbol{E}^{ty/2}\|\|_{\boldsymbol{\mu}}. \tag{4.3}$$

Subsection 4.4 shows that

$$\|\|\boldsymbol{E}^{ty/2}\widehat{\boldsymbol{Q}}\boldsymbol{E}^{ty/2}\|\|_{\boldsymbol{\mu}} \leq \exp\left(t\cdot\boldsymbol{\mu}(y) + \frac{t^2}{2}\cdot\frac{1+\lambda}{1-\lambda}\cdot\frac{(b-a)^2}{4}\right), \tag{4.4}$$

where $\boldsymbol{\mu}(y) \coloneqq (1-\mu)a + \mu b = \pi(f)$. Putting (4.3) and (4.4) together yields

$$\|\|E^{tf/2}\widehat{P}E^{tf/2}\|\|_\pi \leq \exp\left(t\cdot\pi(f) + \frac{t^2}{2}\cdot\frac{1+\lambda}{1-\lambda}\cdot\frac{(b-a)^2}{4}\right). \tag{4.5}$$

Combining (4.5) with (4.1) completes the proof of Theorem 2.1.

León and Perron (2004) and Miasojedow (2014) did not establish (4.1), so their results did not apply to time-dependent functions $f_i$. They did not establish (4.4), so they do not derive explicit bounds for the mgf of $\sum_i f_i(X_i)$ or $\sum_i f(X_i)$.

(Miasojedow, 2014, Lemma 3.9) proves (4.2) in case that $f$ is a simple function taking finitely many possible values. We generalize his result to a general function $f$. Our technique is also different to his: we adopted a simple and clear approach invoking Weyl (1909)'s theorem on essential spectrum, whereas he utilized an analytic argument. These results can be of independent interest.

4.1. *Proof of* (4.1). Lemma 4.2 is devoted to the proof of (4.1). It invokes a few properties of León-Perron operators we proved and summarized in Lemma 4.1. The proof of Lemma 4.1 is deferred to the appendix.

LEMMA 4.1. *Let $\widehat{P}_c = cI + (1-c)\Pi$ with $c \in [0,1)$ be a León-Perron operator on $\mathcal{L}_2(\pi)$. Then the following statement hold.*

*(i) For any bounded function $g$, denote by $M_g$ its multiplication operator. Then*

$$\|g\|_\pi \leq \|\|M_g\widehat{P}_c M_g\|\|_\pi^{1/2}.$$

*Let $P$ be a Markov operator with invariant distribution $\pi$ and absolute spectral gap $1-\lambda$, and let $\widehat{P} = \lambda I + (1-\lambda)\Pi$ be its León-Perron version. The following statements hold.*

(ii) For any $h_1, h_2 \in \mathcal{L}_2(\pi)$,

$$|\langle Ph_1, h_2\rangle_\pi| \leq \langle \widehat{P}h_1, h_1\rangle_\pi^{1/2} \langle \widehat{P}h_2, h_2\rangle_\pi^{1/2}.$$

(iii) For any self-adjoint operators $S_1, S_2$ acting on $\mathcal{L}_2(\pi)$,

$$\||S_1 P S_2\||_\pi \leq \||S_1 \widehat{P} S_1\||_\pi^{1/2} \||S_2 \widehat{P} S_2\||_\pi^{1/2}.$$

LEMMA 4.2. *Let $\{X_i\}_{i\geq 1}$ be a Markov chain with invariant measure $\pi$ and absolute spectral gap $1 - \lambda \in (0, 1]$. Denote by $P$ its Markov operator and by $\widehat{P} = \lambda I + (1-\lambda)\Pi$ the León-Perron version of $P$. Then for any bounded functions $f_i : \mathcal{X} \to [a_i, b_i]$ and any $t \in \mathbb{R}$,*

$$\mathbb{E}_\pi\left[e^{t\sum_{i=1}^n f_i(X_i)}\right] \leq \prod_{i=1}^n \||E^{tf_i/2}\widehat{P}E^{tf_i/2}\||_\pi.$$

PROOF OF LEMMA 4.2. The expression (3.2) implies

$$\mathbb{E}_\pi\left[e^{t\sum_{i=1}^n f_i(X_i)}\right] \leq \|e^{tf_1/2}\|_\pi \prod_{i=1}^{n-1} \||E^{tf_i/2} P E^{tf_{i+1}/2}\||_\pi \|e^{tf_n/2}\|_\pi.$$

Taking $g = e^{tf_1/2}$ or $e^{tf_n/2}$, $S_1 = E^{tf_i/2}$ and $S_2 = E^{tf_{i+1}/2}$ in Lemma 4.1 yields

$$\|e^{tf_1/2}\|_\pi \leq \||E^{tf_1/2}\widehat{P}E^{tf_1/2}\||_\pi^{1/2},$$
$$\|e^{tf_n/2}\|_\pi \leq \||E^{tf_n/2}\widehat{P}E^{tf_n/2}\||_\pi^{1/2},$$
$$\||E^{tf_i/2} P E^{tf_{i+1}/2}\||_\pi \leq \||E^{tf_i/2}\widehat{P}E^{tf_i/2}\||_\pi^{1/2} \||E^{tf_{i+1}/2}\widehat{P}E^{tf_{i+1}/2}\||_\pi^{1/2}.$$

Putting the above displays together completes the proof. $\square$

4.2. *Proof of* (4.2). With Lemma 4.2 in hand, it suffices to consider $\||E^{tf/2}\widehat{P}_c E^{tf/2}\||_\pi$ for function $f : \mathcal{X} \to [a, b]$ and a León-Perron operator $\widehat{P}_c$. Lemma 4.3 shows that $\||E^{tf/2}\widehat{P}_c E^{tf/2}\||_\pi$ can be characterized by the limiting behavior of the mgf of $\sum_{i=1}^n f(\widehat{X}_i)$ when $\{\widehat{X}_i\}_{i\geq 1}$ is driven by $\widehat{P}_c$.

LEMMA 4.3. *Let $\{\widehat{X}_i\}_{i\geq 1}$ be a Markov chain driven by a León-Perron operator $\widehat{P}_c = cI + (1-c)\Pi$ with some $c \in [0, 1)$. Then, for any bounded function $f : \mathcal{X} \to [a, b]$ and any $t \in \mathbb{R}$,*

(i)
$$\mathbb{E}_\pi\left[e^{t\sum_{i=1}^n f(\widehat{X}_i)}\right] \leq \||E^{tf/2}\widehat{P}_c E^{tf/2}\||_\pi^n,$$

(ii)
$$\liminf_{n\to\infty} \frac{1}{n} \log \mathbb{E}_\pi\left[e^{t\sum_{i=1}^n f(\widehat{X}_i)}\right] \geq \log \||E^{tf/2}\widehat{P}_c E^{tf/2}\||_\pi.$$

*Combining the above bounds yields*

$$\lim_{n\to\infty} \frac{1}{n} \log \mathbb{E}_\pi \left[ e^{t \sum_{i=1}^n f(\widehat{X}_i)} \right] = \log \||E^{tf/2} \widehat{P}_c E^{tf/2}\||_\pi.$$

Note that $\widehat{P}_c$ is its own León-Perron operator, Lemma 4.3(i) is a special case of Lemma 4.2 for time-independent functions $f_i = f$ and Markov operator $P = \widehat{P}_c$.

For Lemma 4.3,(ii), we first restrict the attention to the cases in which $f$ is a simple function, namely $f$ takes finitely many possible values, and then extend to a general function by taking the general function as the limit of a sequence of simple functions.

For a simple function, we observe that $E^{tf/2} \widehat{P}_c E^{tf/2}$ behaves like an operator on finite state space, and conclude that the norm of $E^{tf/2} \widehat{P}_c E^{tf/2}$ coincides with its largest eigenvalue. This analysis, summarized in Lemma 4.4, serves as an intermediate step to prove Lemma 4.3(ii).

PROOF OF LEMMA 4.3. (i) Letting $f_i = f$ be time-independent and $P = \widehat{P}_c$ be León-Perron in Lemma 4.2 yields the desired result.

(ii) It is trivial for $t = 0$. For $t \neq 0$, let $\delta(k) = (b-a)/k$ for some $k \in \mathbb{Z}^+$ and define a simple function

$$f_\delta = \begin{cases} a + \delta \lceil (f-a)/\delta \rceil & \text{if } t > 0, \\ a + \delta \lfloor (f-a)/\delta \rfloor & \text{if } t < 0, \end{cases}$$

where $\lfloor \cdot \rfloor$ and $\lceil \cdot \rceil$ are the floor and ceiling function, respectively. Hence $tf_\delta \geq tf \geq tf_\delta - t\delta$.

The self-adjoint operator $\widehat{P}_c$ preserves the non-negativity of $h$ (i.e. $\widehat{P}_c h \geq 0$ if $h \geq 0$), so does self-adjoint operators $E^{tf/2} \widehat{P}_c E^{tf/2}$ and $E^{tf_\delta/2} \widehat{P}_c E^{tf_\delta/2}$. Thus

$$\begin{aligned}
\||E^{tf/2} \widehat{P}_c E^{tf/2}\||_\pi &= \sup_{h:\ \|h\|_\pi = 1} |\langle E^{tf/2} \widehat{P}_c E^{tf/2} h, h \rangle_\pi| \\
&= \sup_{h \geq 0:\ \|h\|_\pi = 1} \langle E^{tf/2} \widehat{P}_c E^{tf/2} h, h \rangle_\pi \\
&\leq \sup_{h \geq 0:\ \|h\|_\pi = 1} \langle E^{tf_\delta/2} \widehat{P}_c E^{tf_\delta/2} h, h \rangle_\pi \quad [tf_\delta \geq tf] \\
&= \sup_{h:\ \|h\|_\pi = 1} |\langle E^{tf_\delta/2} \widehat{P}_c E^{tf_\delta/2} h, h \rangle_\pi| \\
&= \||E^{tf_\delta/2} \widehat{P}_c E^{tf_\delta/2}\||_\pi. \qquad (4.6)
\end{aligned}$$

Suppose we have established (ii) for a simple function like $f_\delta$ in Lemma 4.4

at this moment. Write

$$\liminf_{n\to\infty} \frac{1}{n} \log \mathbb{E}_\pi \left[ e^{t\sum_{i=1}^n f(\widehat{X}_i)} \right]$$
$$\geq \liminf_{n\to\infty} \frac{1}{n} \log \mathbb{E}_\pi \left[ e^{\sum_{i=1}^n tf_\delta(\widehat{X}_i) - n\delta t} \right]$$
$$= \liminf_{n\to\infty} \frac{1}{n} \log \mathbb{E}_\pi \left[ e^{t\sum_{i=1}^n f_\delta(\widehat{X}_i)} \right] - \delta t$$
$$\geq \log \||E^{tf_\delta/2} \widehat{P}_c E^{tf_\delta/2}\||_\pi - \delta t \qquad \text{[Lemma 4.4(iv)]}$$
$$\geq \log \||E^{tf/2} \widehat{P}_c E^{tf/2}\||_\pi - \delta t. \qquad \text{[By (4.6)]}$$

Letting $\delta$ tend to 0 ($k$ tending to $\infty$) completes the proof.

$\square$

Now it remains to prove Lemma 4.3(ii) for a simple function. This is completed in Lemma 4.4.

LEMMA 4.4. *Let $\widehat{P}_c = cI + (1-c)\Pi$ be a León-Perron operator. Let $f$ be a simple function on $\mathcal{X}$ $\pi$-a.e, that is, there exists a finite set $\{\beta_1, \ldots, \beta_k\}$ with $\beta_1 > \cdots > \beta_k$ such that $f^{-1}(\beta_j) := \{x \in \mathcal{X} : f(x) = \beta_j\}$ satisfies*

$$\pi(f^{-1}(\beta_j)) > 0, \ \forall 1 \leq j \leq k; \qquad \sum_{j=1}^k \pi(f^{-1}(\beta_j)) = 1.$$

*Let*

$$F(r) = \pi\left(\frac{(1-c)e^f}{r - ce^f}\right) = \sum_{j=1}^k \frac{(1-c)e^{\beta_j}}{r - ce^{\beta_j}} \pi(f^{-1}(\beta_j)).$$

*The following statements hold.*

(i) *Each solution $r_\star$ to $F(r_\star) = 1$ is an eigenvalue of $E^{f/2} \widehat{P}_c E^{f/2}$ with eigenfunction*

$$h_\star = \frac{(1-c)e^{f/2}}{r_\star - ce^f}.$$

*Let $\beta_0 = \infty$. There are $k$ such solutions $r_j \in (ce^{\beta_j}, ce^{\beta_{j-1}})$ for $j = 1, \ldots, k$.*

(ii) *Denote by $\mathcal{S}(E^{f/2} \widehat{P}_c E^{f/2})$ the spectrum of operator $E^{f/2} \widehat{P}_c E^{f/2}$. Then*

$$\{r_j : 1 \leq j \leq k\} \subseteq \mathcal{S}(E^{f/2} \widehat{P}_c E^{f/2})$$
$$\subseteq \{ce^{\beta_j} : 1 \leq j \leq k\} \cup \{r_j : 1 \leq j \leq k\}.$$

(iii) *Recall that $r_1$ is the largest eigenvalue of $E^{f/2} \widehat{P}_c E^{f/2}$. Then*

$$\||E^{f/2} \widehat{P}_c E^{f/2}\||_\pi = r_1.$$

(iv) *Let $\{\widehat{X}_i\}_{i\geq 1}$ be a Markov chain driven by $\widehat{P}_c$ then*

$$\liminf_{n\to\infty} \frac{1}{n} \log \mathbb{E}_\pi \left[ e^{\sum_{i=1}^n f(\widehat{X}_i)} \right] \geq \log r_1 = \log \||E^{f/2} \widehat{P}_c E^{f/2}\||_\pi.$$

PROOF OF LEMMA 4.4. (i) Note that $\pi(e^{f/2}h_\star) = F(r_\star) = 1$. We have

$$
\begin{aligned}
&E^{f/2}\widehat{P}_c E^{f/2} h_\star(x) - r_\star h_\star(x) \\
&= ce^{f(x)}h_\star(x) + (1-c)e^{f(x)}\pi(e^{f/2}h_\star) - r_\star h_\star(x) \\
&= ce^{f(x)}h_\star(x) + (1-c)e^{f(x)}F(r_\star) - r_\star h_\star(x).
\end{aligned}
$$

Plugging $F(r_\star) = 1$ and $h_\star = (1-c)e^{f/2}/(r_\star - ce^f)$ into the last line yields 0. Thus $r_\star$ is an eigenvalue of $E^{f/2}\widehat{P}_c E^{f/2}$ with eigenfunction $h_\star$.

On each interval $(ce^{\beta_j}, ce^{\beta_{j-1}})$, function $F(r)$ decreases to $-\infty$ (or 0 if $j=1$) as $r \uparrow ce^{\beta_{j-1}}$, and increases to $+\infty$ as $r \downarrow ce^{\beta_j}$. Thus there exists $r_j \in (ce^{\beta_j}, ce^{\beta_{j-1}})$ such that $F(r_j) = 1$.

(ii) The operator $E^{f/2}\widehat{P}_c E^{f/2}$ is self-adjoint, thus its spectrum consists of the discrete spectrum $\mathcal{S}_d(E^{f/2}\widehat{P}_c E^{f/2})$, which are isolated eigenvalues of finite multiplicity, and the essential spectrum $\mathcal{S}_{ess}(E^{f/2}\widehat{P}_c E^{f/2})$.

We first prove $\mathcal{S}_d(E^{f/2}\widehat{P}_c E^{f/2}) \subseteq \{ce^{\beta_j} : 1 \le j \le k\} \cup \{r_j : 1 \le j \le k\}$ by showing that any eigenvalue $r_\star$ belongs to either $\{ce^{\beta_j} : 1 \le j \le k\}$ or $\{r_j : 1 \le j \le k\}$. Consider any pair of eigenvalue $r_\star$ and non-zero eigenfunction $h_\star$ such that

$$
\begin{aligned}
r_\star h_\star(x) &= E^{f/2}\widehat{P}_c E^{f/2} h_\star(x) \\
&= ce^{f(x)}h_\star(x) + (1-c)e^{f(x)/2}\pi(e^{f/2}h_\star), \quad \pi\text{-a.e. } x.
\end{aligned}
$$

There are two possibilities. If $\pi(e^{f/2}h_\star) = 0$ then the above display implies that $(r_\star - ce^{f(x)})h_\star(x) = 0$ for $\pi$-a.e. $x$. There exists at least one $j$ such that $h_\star(x)$ is not identically zero on the set $f^{-1}(\beta_j)$, implying $r_\star = ce^{\beta_j}$. If $\pi(e^{f/2}h_\star) \ne 0$ and $r_\star \notin \{ce^{\beta_j} : 1 \le j \le k\}$ then we have

$$
h_\star(x) = \frac{(1-c)e^{f(x)/2}}{r_\star - ce^{f(x)}}\pi(e^{f/2}h_\star), \quad \pi\text{-a.e. } x.
$$

Multiplying both sides with $e^{f(x)/2}$, taking expectation of both sides with respect to $\pi$, and dividing both sides by $\pi(e^{f/2}h_\star)$ yields $F(r_\star) = 1$, that means $r_\star \in \{r_j : 1 \le j \le k\}$.

We next show that $\mathcal{S}_{ess}(E^{f/2}\widehat{P}_c E^{f/2}) \subseteq \{ce^{\beta_j} : 1 \le j \le k\}$ by Weyl (1909)'s theorem on essential spectrum. Write

$$
\widehat{P}_c = cE^f + (1-c)E^{f/2}\Pi E^{f/2}
$$

in the form of a self-adjoint operator $cE^f$ perturbed by another self-adjoint operator $(1-c)E^{f/2}\Pi E^{f/2}$. The perturbation operator $(1-c)E^{f/2}\Pi E^{f/2}$ is of finite rank, as

$$
(1-c)E^{f/2}\Pi E^{f/2} h = (1-c)e^{f/2}\pi(e^{f/2}h) \in \{\xi e^{f/2} : \xi \in \mathbb{R}\},
$$

and thus compact. Weyl (1909) asserts that the essential spectrum of a self-adjoint operator is invariant to the perturbation of a self-adjoint, compact operator. Hence $\widehat{P}_c$ shares the same essential spectrum with $cE^f$. Recall that $cE^f$ is the multiplication operator of function $ce^f$. Its spectrum is the essential range of $ce^f$, which is $\{ce^{\beta_j} : 1 \leq j \leq k\}$. Thus

$$\mathcal{S}_{\mathrm{ess}}(E^{f/2}\widehat{P}_c E^{f/2}) = \mathcal{S}_{\mathrm{ess}}(cE^f) \subseteq \{ce^{\beta_j} : 1 \leq j \leq k\}.$$

Combining (i) and results for $\mathcal{S}_{\mathrm{d}}(E^{f/2}\widehat{P}_c E^{f/2})$ and $\mathcal{S}_{\mathrm{ess}}(E^{f/2}\widehat{P}_c E^{f/2})$ completes the proof.

(iii) By (ii), $r_1$ is the spectral radius of $E^{f/2}\widehat{P}_c E^{f/2}$. Recall that $E^{f/2}\widehat{P}_c E^{f/2}$ is self-adjoint. Thus $r_1 = |\!|\!| E^{f/2}\widehat{P}_c E^{f/2} |\!|\!|_\pi$.

(iv) By (i), eigenvalue $r_1$ associates with eigenfunction

$$h_1 = \frac{(1-c)e^{f/2}}{r_1 - ce^f}$$

and $\langle h_1, e^{f/2}\rangle_\pi = F(r_1) = 1$. Let $\widetilde{h}_1$ be the projection of $e^{f/2}$ onto $h_1$. It is elementary that

$$\widetilde{h}_1 := \left\langle \frac{h_1}{\|h_1\|_\pi}, e^{f/2}\right\rangle_\pi \frac{h_1}{\|h_1\|_\pi} = \frac{h_1}{\|h_1\|_\pi^2},$$

and

$$\langle e^{f/2} - \widetilde{h}_1, \widetilde{h}_1\rangle_\pi = 0.$$

$h_1$ is the eigenfunction of the self-adjoint operator $E^{f/2}\widehat{P}_c E^{f/2}$, thus

$$\begin{aligned}0 &= \langle e^{f/2} - \widetilde{h}_1, (E^{f/2}\widehat{P}_c E^{f/2})^{n-1}\widetilde{h}_1\rangle_\pi \\ &= \langle (E^{f/2}\widehat{P}_c E^{f/2})^{n-1}(e^{f/2} - \widetilde{h}_1), \widetilde{h}_1\rangle_\pi.\end{aligned}$$

The self-adjoint operator $E^{f/2}\widehat{P}_c E^{f/2}$ is positive semi-definite, so

$$\left\langle e^{f/2} - \widetilde{h}_1, (E^{f/2}\widehat{P}_c E^{f/2})^{n-1}(e^{f/2} - \widetilde{h}_1)\right\rangle_\pi \geq 0.$$

Thus

$$\begin{aligned}\mathbb{E}_\pi\left[e^{\sum_{i=1}^n f(\widehat{X}_i)}\right] &= \left\langle e^{f/2}, (E^{f/2}\widehat{P}_c E^{f/2})^{n-1} e^{f/2}\right\rangle_\pi \\ &\geq \left\langle \widetilde{h}_1, (E^{f/2}\widehat{P}_c E^{f/2})^{n-1}\widetilde{h}_1\right\rangle_\pi \\ &= r_1^{n-1}/\|h_1\|_\pi^2,\end{aligned}$$

implying the desired result.

$\square$

4.3. *Proof of* (4.3). This subsection is devoted to the proof of (4.3). First, we have Lemma 4.5 as a corollary of Lemma 4.3 for the two-state chain $\{\widehat{Y}_i\}_{i\geq 1}$.

LEMMA 4.5 (Application of Lemma 4.3 to the two-state chain). *Let $\{\widehat{Y}_i\}_{i\geq 1}$ be the two-state Markov chain in Definition 4.1. Recall that $\boldsymbol{E}^{ty/2} = diag(e^{ta/2}, e^{tb/2})$ and that $\|\|\boldsymbol{T}\|\|_{\boldsymbol{\mu}}$ denotes the operator norm induced by the $\boldsymbol{\mu}$-weighted vector norm for a $2 \times 2$ matrix $\boldsymbol{T}$. Then*

$$\lim_{n\to\infty} \frac{1}{n} \log \mathbb{E}_{\boldsymbol{\mu}} \left[ e^{t \sum_{i=1}^n \widehat{Y}_i} \right] = \log \|\|\boldsymbol{E}^{ty/2} \widehat{\boldsymbol{Q}} \boldsymbol{E}^{ty/2}\|\|_{\boldsymbol{\mu}}.$$

PROOF OF LEMMA 4.5. We abuse the bolded greek letter $\boldsymbol{\mu}$ to denote the invariant measure of the two-state Markov chain. Recall that transition probability matrix $\widehat{\boldsymbol{Q}}$ is its own León-Perron version on the Hilbert space $\mathcal{L}_2(\boldsymbol{\mu}) = \mathcal{L}_2(\mathcal{Y}, 2^{\mathcal{Y}}, \boldsymbol{\mu})$, where

$$\mathcal{Y} = \{a, b\}, \quad 2^{\mathcal{Y}} = \{\emptyset, \{a\}, \{b\}, \{a, b\}\}, \quad \boldsymbol{\mu}(\{b\}) = \mu.$$

This lemma is merely an application of Lemma 4.3 to the two-state Markov chain $\{\widehat{Y}_i\}_{i\geq 1}$ on the Hilbert space $\mathcal{L}_2(\boldsymbol{\mu})$ and $f$ in Lemma 4.3 replaced by the identity function $y \in \mathcal{Y} \mapsto y$. □

Lemma 4.6, taken from León and Perron (2004), asserts that the two-state Markov chain $\{\widehat{Y}_i\}_{i\geq 1}$ is the extremal case of Markov chains with absolute spectral gap $1 - \lambda$. Putting Lemmas 4.3, 4.5 and 4.6 together concludes (4.3).

LEMMA 4.6 (Theorem 2 in León and Perron (2004)). *Let $\{\widehat{X}_i\}_{i\geq 1}$ be a Markov chain driven by the León-Perron operator $\widehat{P} = \lambda I + (1-\lambda)\Pi$ with some $\lambda \in [0, 1)$. For a bounded function $f : \mathcal{X} \to [a, b]$, let $\{\widehat{Y}_i\}_{i\geq 1}$ be the two-state Markov chain in Definition 4.1. If $\mu = (\pi(f) - a)/(b - a)$ then, for any convex function $\Psi : \mathbb{R} \to \mathbb{R}$,*

$$\mathbb{E}_\pi \left[ \Psi \left( \sum_{i=1}^n f(\widehat{X}_i) \right) \right] \leq \mathbb{E}_{\boldsymbol{\mu}} \left[ \Psi \left( \sum_{i=1}^n \widehat{Y}_i \right) \right].$$

*In particular, for $\Psi : z \mapsto \exp(tz)$,*

$$\mathbb{E}_\pi \left[ e^{t \sum_{i=1}^n f(\widehat{X}_i)} \right] \leq \mathbb{E}_{\boldsymbol{\mu}} \left[ e^{t \sum_{i=1}^n \widehat{Y}_i} \right].$$

4.4. *Proof of* (4.4). So far the task of bounding $\|\|E^{tf/2} \widehat{P} E^{tf/2}\|\|_\pi$ has been reduced to that of bounding $\|\|\boldsymbol{E}^{ty/2} \widehat{\boldsymbol{Q}} \boldsymbol{E}^{ty/2}\|\|_{\boldsymbol{\mu}}$ in the two-state Markov chain system on the state space $\{a, b\}$. The latter is attacked by Lemma 4.7.

LEMMA 4.7. *Use the same notations in Lemma 4.5.*

*(i) Let $\theta(t)$ be the largest eigenvalue of matrix $\boldsymbol{E}^{ty/2} \widehat{\boldsymbol{Q}} \boldsymbol{E}^{ty/2}$, then*

$$\|\|\boldsymbol{E}^{ty/2} \widehat{\boldsymbol{Q}} \boldsymbol{E}^{ty/2}\|\|_{\boldsymbol{\mu}} = \theta(t).$$

*(ii) Let $\boldsymbol{\mu}(y) := (1-\mu)a + \mu b$ and recall that $\alpha : \lambda \mapsto (1+\lambda)/(1-\lambda)$. Then*

$$\theta(t) \leq \widetilde{\theta}(t) := \exp\left( t \cdot \boldsymbol{\mu}(y) + \frac{t^2}{2} \cdot \alpha(\lambda) \cdot \frac{(b-a)^2}{4} \right).$$

PROOF OF LEMMA 4.7. (i) It follows from Lemma 4.4(iii), as function $y \in \mathcal{Y} \mapsto y$ takes 2 possible values $a$ and $b$. This is also straightforward by Frobenius-Perron theorem.

(ii) Let
$$p = \frac{\lambda + (1-\lambda)\mu}{1+\lambda}, \quad 1-p = \frac{\lambda + (1-\lambda)(1-\mu)}{1+\lambda}.$$

The largest eigenvalue $\theta(t)$ is the right solution to the following quadratic equation
$$\begin{aligned} 0 &= \det(\theta \boldsymbol{I} - \boldsymbol{E}^{ty/2}\widehat{\boldsymbol{Q}}\boldsymbol{E}^{ty/2}) \\ &= \theta^2 - \left[(\lambda + (1-\lambda)(1-\mu))e^{ta} + (\lambda + (1-\lambda)\mu)e^{tb}\right]\theta + \lambda e^{ta+tb} \\ &= \theta^2 - (1+\lambda)\left[(1-p)e^{ta} + pe^{tb}\right]\theta + \lambda e^{ta+tb}.\end{aligned}$$

It suffices to show $\widetilde{\theta}(t) = \exp\left(t \cdot \boldsymbol{\mu}(y) + \frac{t^2}{2} \cdot \alpha(\lambda) \cdot \frac{(b-a)^2}{4}\right)$ satisfies

$$\widetilde{\theta}(t)^2 - (1+\lambda)\left[(1-p)e^{ta} + pe^{tb}\right]\widetilde{\theta}(t) + \lambda e^{ta+tb} \geq 0, \text{ and} \qquad (4.7)$$
$$\widetilde{\theta}(t)^2 \geq \lambda e^{ta+tb}. \qquad (4.8)$$

An equivalent form of (4.7) is as follows.
$$\frac{\widetilde{\theta}(t) + \lambda e^{ta+tb}\widetilde{\theta}(t)^{-1}}{1+\lambda} \geq (1-p)e^{ta} + pe^{tb}. \qquad (4.9)$$

Using the convexity of function $z \mapsto e^z$, the left hand side of (4.9) is lower bounded as

$$\begin{aligned} &\frac{\widetilde{\theta}(t) + \lambda e^{ta+tb}\widetilde{\theta}(t)^{-1}}{1+\lambda} \\ &= \frac{\exp(t\boldsymbol{\mu}(y) + \alpha(\lambda)(b-a)^2 t^2/8) + \lambda \exp(at + bt - t\boldsymbol{\mu}(y) - \alpha(\lambda)(b-a)^2 t^2/8)}{1+\lambda} \\ &\geq \exp\left(\frac{t\boldsymbol{\mu}(y) + \alpha(\lambda)(b-a)^2 t^2/8 + \lambda at + \lambda bt - \lambda t\boldsymbol{\mu}(y) - \lambda \alpha(\lambda)(b-a)^2 t^2/8}{1+\lambda}\right) \\ &= \exp\left(t \cdot \frac{(1-\lambda)\boldsymbol{\mu}(y) + \lambda a + \lambda b}{1+\lambda} + \frac{(b-a)^2 t^2}{8} \cdot \frac{(1-\lambda)\alpha(\lambda)}{1+\lambda}\right) \\ &= \exp\left(t \cdot [(1-p)a + pb] + \frac{(b-a)^2 t^2}{8}\right).\end{aligned}$$

The right hand side of (4.9) is the mgf of a Bernoulli random variable $Z$ with $\mathbb{P}(Z = a) = 1-p$ and $\mathbb{P}(Z = b) = p$, which is upper bounded by the classical Hoeffding's lemma as
$$(1-p)e^{ta} + pe^{tb} \leq \exp\left(t \cdot [(1-p)a + pb] + \frac{(b-a)^2 t^2}{8}\right).$$

Comparing both sides yields (4.9) and (4.7). On the other hand, (4.8) holds as

$$\log\left(\widetilde{\theta}(t)^2 e^{-ta-tb}\right) = \frac{\alpha(\lambda)(b-a)^2 t^2}{4} + (2\mu-1)(b-a)t$$
$$\geq -\frac{(2\mu-1)^2}{\alpha(\lambda)}$$
$$\geq -\frac{1}{\alpha(\lambda)} = -\frac{1-\lambda}{1+\lambda} \geq \log \lambda.$$

□

**5. Proof of Theorem 2.2.** The proof uses properties of *numerical range* of the additively-reversiblized operator $R = (P + P^*)/2$.

DEFINITION 5.1 (Numerical range and radius). *The numerical range of an operator $T$ on a Hilbert space $\mathcal{H}$ is defined as*

$$\mathcal{W}(T|\mathcal{H}) := \{\langle Th, h\rangle : h \in \mathcal{H}, \|h\| = 1\}.$$

*The numerical radius of $\mathcal{W}(T|\mathcal{H})$ is defined as*

$$\rho[\mathcal{W}(T|\mathcal{H})] := \sup\{|w| : w \in \mathcal{W}(T|\mathcal{H})\}.$$

LEMMA 5.1 (Gustafson and Rao (1997), p. 39). *If a bounded normal operator $T_1$ on a Hilbert space $\mathcal{H}$ commutes with $T_2$ then*

$$\rho[\mathcal{W}(T_1 T_2|\mathcal{H})] \leq \rho[\mathcal{W}(T_1|\mathcal{H})]\rho[\mathcal{W}(T_2|\mathcal{H})].$$

*In particular, for any self-adjoint operator $S$,*

$$\rho[\mathcal{W}(S^n|\mathcal{H})] \leq (\rho[\mathcal{W}(S|\mathcal{H})])^n.$$

Since $R = (P + P^*)/2$ is self-adjoint, we have

$$\lambda_{\rm r}(R) = \sup\{w : w \in \mathcal{W}(R|\mathcal{L}_2^0(\pi))\}.$$

With these results in hand, we proceed to prove Theorem 2.2.

PROOF OF THEOREM 2.2. Without loss of generality, we assume $\pi(f_i) = 0$. From (3.2), Lemma 5.1, and the fact that operator $E^{tf/2} P E^{tf/2}$ preserves non-negativity, it follows that

$$\mathbb{E}_\pi\left[e^{t\sum_{i=1}^n f(X_i)}\right] = \langle e^{tf/2}, (E^{tf/2} P E^{tf/2})^{n-1} e^{tf/2}\rangle_\pi$$
$$\leq \sup_{\|h\|_\pi=1} |\langle (E^{tf/2} P E^{tf/2})^{n-1} h, h\rangle| \|e^{tf/2}\|_2^2$$
$$\leq \left(\sup_{\|h\|_\pi=1} |\langle E^{tf/2} P E^{tf/2} h, h\rangle|\right)^{n-1} \|e^{tf/2}\|_2^2$$
$$= \left(\sup_{\|h\|_\pi=1} \langle E^{tf/2} P E^{tf/2} h, h\rangle\right)^{n-1} \|e^{tf/2}\|_2^2. \qquad (5.1)$$

It is elementary that
$$\langle E^{tf/2} P E^{tf/2} h, h \rangle = \langle E^{tf/2} R E^{tf/2} h, h \rangle. \tag{5.2}$$

Let $\widehat{R} = (\lambda_r \vee 0) I + (1 - (\lambda_r \vee 0)) \Pi$ be a León-Perron operator. From $R = (I - \Pi) R (I - \Pi) + \Pi$ (since $R\Pi = \Pi R = \Pi$) and self-adjointness of $E^{tf/2}$ and $I - \Pi$ it follows that

$$\begin{aligned}
\langle E^{tf/2} R E^{tf/2} h, h \rangle_\pi &= \langle E^{tf/2} (I - \Pi) h, R(I - \Pi) E^{tf/2} h \rangle_\pi + \pi(e^{tf/2} h)^2 \\
&\leq \lambda_r \| (I - \Pi) E^{tf/2} h \|_\pi^2 + \pi(E^{tf/2} h)^2 \\
&\leq (\lambda_r \vee 0) \| (I - \Pi) E^{tf/2} h \|_\pi^2 + \pi(E^{tf/2} h)^2 \\
&= \langle E^{tf/2} \widehat{R} E^{tf/2} h, h \rangle_\pi, \quad \forall h \in \mathcal{L}_2(\pi),
\end{aligned}$$

implying

$$\sup_{\|h\|_\pi = 1} \langle E^{tf/2} R E^{tf/2} h, h \rangle_\pi \leq \sup_{\|h\|_\pi = 1} \langle E^{tf/2} \widehat{R} E^{tf/2} h, h \rangle_\pi$$
$$= \| E^{tf/2} \widehat{R} E^{tf/2} \|_\pi. \tag{5.3}$$

By Lemma 4.1(i),
$$\| e^{tf/2} \|_2^2 \leq \| E^{tf/2} \widehat{R} E^{tf/2} \|_\pi. \tag{5.4}$$

Putting (5.1), (5.2), (5.3) and (5.4) together yields

$$\mathbb{E}_\pi \left[ e^{t \sum_{i=1}^n f(X_i)} \right] \leq \| E^{tf/2} \widehat{R} E^{tf/2} \|_\pi^n.$$

Noting that (4.5) holds for a general León-Perron operator, we have

$$\| E^{tf/2} \widehat{R} E^{tf/2} \|_\pi \leq \exp\left( \frac{t^2}{2} \cdot \alpha(\lambda_r \vee 0) \cdot \frac{(b-a)^2}{4} \right),$$

concluding the proof. $\square$

**6. Applications.** This section applies our theorems to streamline the analyses of six statistics and machine learning problems involving Markov-dependence. Let us collect more notations for vector and matrix norms. For a vector $\boldsymbol{z}$, let $\|\boldsymbol{z}\|_1$ and $\|\boldsymbol{z}\|$ denote the $\ell_1$-norm and $\ell_2$-norm, respectively. For a matrix $\boldsymbol{S}$, let $\|\boldsymbol{S}\|_1$ and $\|\boldsymbol{S}\|$ denote the operator norms induced by the $\ell_1$-norm and $\ell_2$-norm of vectors. If $\boldsymbol{S}$ is a $d \times d$ symmetric matrix then it is elementary that $\|\boldsymbol{S}\| \leq \|\boldsymbol{S}\|_1 \leq d \max_{ij} |\boldsymbol{S}_{ij}|$.

**6.1. Linear regression.** We consider a linear regression model with Markov-dependent samples.

$$y_i = \boldsymbol{f}(X_i)' \boldsymbol{\beta}_\star + \varepsilon_i, \quad \text{for } i = 1, \ldots, n, \tag{6.1}$$

where $X_i$'s are samples from a stationary Markov chain $\{X_i\}_{i \geq 1}$ with invariant measure $\pi$ and right spectral gap $1 - \lambda_r$; the vector of noise terms $\boldsymbol{\varepsilon} = (\varepsilon_1, \ldots, \varepsilon_n)'$ is sub-Gaussian with variance proxy $\sigma^2$ and independent from $\{X_i\}_{i=1}^n$; and, $\boldsymbol{f} = (f_1, \ldots, f_d)'$ is a known $d$-dimensional bounded feature mapping. Without loss of generality, we assume $\sup_{x \in \mathcal{X}} |f_j(x)| \leq 1$ for all $j = 1, \ldots, d$. Let $\boldsymbol{F} = [\boldsymbol{f}(X_1), \ldots, \boldsymbol{f}(X_n)]'$ be the $n \times d$ random design matrix, and let $\boldsymbol{F}_j = [f_j(X_1), \ldots, f_j(X_n)]'$ be the $j$-th column of $\boldsymbol{F}$. Also, let

$$\boldsymbol{\Sigma} = \int \boldsymbol{f}(x)\boldsymbol{f}(x)' \pi(dx) = [\pi(f_j f_k)]_{1 \leq j,k \leq d},$$

$$\widehat{\boldsymbol{\Sigma}} = \frac{1}{n} \sum_{i=1}^n \boldsymbol{f}(X_i)\boldsymbol{f}(X_i)' = \left[\frac{1}{n} \sum_{i=1}^n f_j(X_i) f_k(X_i)\right]_{1 \leq j,k \leq d}, \quad (6.2)$$

$$\epsilon_n = \max_{j,k} \left|\widehat{\boldsymbol{\Sigma}}_{j,k} - \boldsymbol{\Sigma}_{j,k}\right|.$$

We attack the low-dimensional case with $d < n$ in this subsection, and leaves the high-dimensional case to the next subsection. In the low-dimensional case, the ordinary least squared estimator is given by

$$\widehat{\boldsymbol{\beta}} = \widehat{\boldsymbol{\Sigma}}^{-1} \left(\boldsymbol{F}'\boldsymbol{y}/n\right), \text{ and } \widehat{\boldsymbol{\beta}} - \boldsymbol{\beta}_\star = \widehat{\boldsymbol{\Sigma}}^{-1}(\boldsymbol{F}'\boldsymbol{\varepsilon}/n).$$

We would like to derive a non-asymptotic error bound for $\widehat{\boldsymbol{\beta}}$ in $\ell_2$-norm. Let us first consider $\boldsymbol{F}'\boldsymbol{\varepsilon}/n$. To this end, $\boldsymbol{F}_j'\boldsymbol{\varepsilon}$ is sub-Gaussian with variance proxy $n\sigma^2$, since $\sup_{x \in \mathcal{X}} |f_j(x)| \leq 1$ and $\boldsymbol{\varepsilon}$ is sub-Gaussian with variance proxy $\sigma^2$. It follows that

$$\mathbb{P}\left(\|\boldsymbol{F}'\boldsymbol{\varepsilon}/n\| \geq \epsilon\right) \leq \sum_{j=1}^d \mathbb{P}\left(|\boldsymbol{F}_j'\boldsymbol{\varepsilon}| \geq n\epsilon/\sqrt{d}\right) \leq 2d \exp\left(-\frac{n\epsilon^2}{2d\sigma^2}\right).$$

The remaining analyses need two technical lemmas as follows. The proof of Lemma 6.1 is concluded by simply applying Theorem 2.2 to each of $d^2$ elements of $\widehat{\boldsymbol{\Sigma}} - \boldsymbol{\Sigma}$ and a union bound over all $d^2$ elements. Lemma 6.2 is due to Sun, Tan, Liu and Zhang (2017).

LEMMA 6.1. *Let $\{X_i\}_{i \geq 1}$ be a stationary Markov chain with invariant measure $\pi$ and right spectral gap $1 - \lambda_r$. Let $\boldsymbol{f} = (f_1, \ldots, f_d)'$ be a $d$-dimensional feature mapping with $\sup_{x \in \mathcal{X}} |f_j(x)| \leq 1$ for all $j = 1, \ldots, d$. Refer notations $\boldsymbol{\Sigma}$, $\widehat{\boldsymbol{\Sigma}}$ and $\epsilon_n$ to (6.2). Then*

$$\mathbb{P}(\epsilon_n \geq \epsilon) \leq 2d^2 \exp\left(-\frac{n\epsilon^2}{2\alpha(\lambda_r \vee 0)}\right).$$

*From the fact that $\|\widehat{\boldsymbol{\Sigma}} - \boldsymbol{\Sigma}\| \leq d\epsilon_n$, it follows that*

$$\mathbb{P}\left(\|\widehat{\boldsymbol{\Sigma}} - \boldsymbol{\Sigma}\| \geq \epsilon\right) \leq 2d^2 \exp\left(-\frac{n\epsilon^2}{2d^2\alpha(\lambda_r \vee 0)}\right).$$

LEMMA 6.2 (Lemma E.4 in Sun, Tan, Liu and Zhang (2017)). *If $\boldsymbol{\Sigma}_1, \boldsymbol{\Sigma}_2 \in \mathbb{R}^{d \times d}$ are invertible and $\|\boldsymbol{\Sigma}_1^{-1}\| \|\boldsymbol{\Sigma}_1 - \boldsymbol{\Sigma}_2\| < 1$ then*

$$\|\boldsymbol{\Sigma}_1^{-1} - \boldsymbol{\Sigma}_2^{-1}\| \leq \frac{\|\boldsymbol{\Sigma}_1^{-1}\|^2 \|\boldsymbol{\Sigma}_1 - \boldsymbol{\Sigma}_2\|}{1 - \|\boldsymbol{\Sigma}_1^{-1}\| \|\boldsymbol{\Sigma}_1 - \boldsymbol{\Sigma}_2\|}.$$

With these pieces in hand, we are ready to bound $\|\widehat{\boldsymbol{\beta}}-\boldsymbol{\beta}_\star\|$. Given that $\|\|\boldsymbol{\Sigma}^{-1}\|\|\|\|\widehat{\boldsymbol{\Sigma}} - \boldsymbol{\Sigma}\|\| \leq 1 - \eta \in (0,1)$, it follows from Lemma 6.2 that

$$\|\widehat{\boldsymbol{\beta}} - \boldsymbol{\beta}_\star\| \leq \|\|\widehat{\boldsymbol{\Sigma}}^{-1} - \boldsymbol{\Sigma}^{-1}\|\| \|\boldsymbol{F}'\boldsymbol{\varepsilon}/n\| + \|\|\boldsymbol{\Sigma}^{-1}\|\| \|\boldsymbol{F}'\boldsymbol{\varepsilon}/n\|$$

$$\leq \frac{\|\|\boldsymbol{\Sigma}^{-1}\|\|^2 \|\|\widehat{\boldsymbol{\Sigma}} - \boldsymbol{\Sigma}\|\|}{1 - \|\|\boldsymbol{\Sigma}^{-1}\|\| \|\|\widehat{\boldsymbol{\Sigma}} - \boldsymbol{\Sigma}\|\|} \|\boldsymbol{F}'\boldsymbol{\varepsilon}/n\| + \|\|\boldsymbol{\Sigma}^{-1}\|\| \|\boldsymbol{F}'\boldsymbol{\varepsilon}/n\|$$

$$\leq \frac{\|\|\boldsymbol{\Sigma}^{-1}\|\|}{1 - \|\|\boldsymbol{\Sigma}^{-1}\|\| \|\|\widehat{\boldsymbol{\Sigma}} - \boldsymbol{\Sigma}\|\|} \|\boldsymbol{F}'\boldsymbol{\varepsilon}/n\|$$

$$\leq \eta^{-1} \|\|\boldsymbol{\Sigma}^{-1}\|\| \|\boldsymbol{F}'\boldsymbol{\varepsilon}/n\|.$$

Thus, when $\eta = 1/2$,

$$\mathbb{P}\left(\|\widehat{\boldsymbol{\beta}} - \boldsymbol{\beta}_\star\| \geq \epsilon\right) \leq \mathbb{P}\left(\|\|\boldsymbol{\Sigma}^{-1}\|\| \|\|\widehat{\boldsymbol{\Sigma}} - \boldsymbol{\Sigma}\|\| \geq 1/2\right) + \mathbb{P}\left(\|\|\boldsymbol{\Sigma}^{-1}\|\| \|\boldsymbol{F}'\boldsymbol{\varepsilon}/n\| \geq \epsilon/2\right)$$

$$\leq 2d^2 \exp\left(-\frac{n\alpha(\lambda_\mathrm{r} \vee 0)^{-1}}{8d^2 \|\|\boldsymbol{\Sigma}^{-1}\|\|^2}\right) + 2d \exp\left(-\frac{n\epsilon^2}{8d\sigma^2 \|\|\boldsymbol{\Sigma}^{-1}\|\|^2}\right).$$

With probability at least $1 - 4\delta$,

$$\|\widehat{\boldsymbol{\beta}} - \boldsymbol{\beta}_\star\| \leq \sigma \|\|\boldsymbol{\Sigma}^{-1}\|\| \sqrt{8d(\log(1/\delta) + \log d)/n},$$

if

$$n \geq 8\alpha(\lambda_\mathrm{r} \vee 0)d^2 \|\|\boldsymbol{\Sigma}^{-1}\|\|^2 (\log(1/\delta) + 2\log d).$$

In contrast, the case of the i.i.d. data samples requires

$$n \geq 8d^2 \|\|\boldsymbol{\Sigma}^{-1}\|\|^2 (\log(1/\delta) + 2\log d).$$

6.2. *Sparse regression.* Proceed to the $s$-sparse instance of the linear model (6.1) in the high-dimensional regime with $d \gg n$. Here we assume $\boldsymbol{\beta}_\star$ is $s$-sparse, i.e. its support, denoted by $\mathcal{S}$, is of size $|\mathcal{S}| = s$. The lasso solution with a regularization penalty $w_n$ is given by

$$\boldsymbol{\beta}_{w_n} \in \arg\min_{\boldsymbol{\beta} \in \mathbb{R}^d} \left\{ \frac{1}{2n} \|\boldsymbol{y} - \boldsymbol{F}\boldsymbol{\beta}\|^2 + w_n \|\boldsymbol{\beta}\|_1 \right\}.$$

See Tibshirani (1996); Fan and Li (2001); Zou (2006); Zou and Li (2008); Fan and Lv (2011) for a perspective on this model.

Corollary 2 in Negahban et al. (2012) asserts that, with probability at least $1 - c_1 \exp\left(-c_2 n w_n^2\right)$, any optimal solution $\boldsymbol{\beta}_{w_n}$ with $w_n = 4\sigma \sqrt{\log d/n}$ satisfies the bounds

$$\|\boldsymbol{\beta}_{w_n} - \boldsymbol{\beta}_\star\| \leq \frac{8\sigma}{\kappa} \sqrt{s \log d/n}, \quad \|\boldsymbol{\beta}_{w_n} - \boldsymbol{\beta}_\star\|_1 \leq \frac{24\sigma}{\kappa} s \sqrt{\log d/n},$$

if the following two conditions hold. One is the *normalized column* condition: each $j$-th column $\boldsymbol{F}_j$ of the design matrix satisfies

$$\|\boldsymbol{F}_j\| \leq \sqrt{n}.$$

The other is the *restricted eigenvalue* condition: there exists some absolute constant $\kappa$ such that for all $\boldsymbol{\beta}$ satisfying $\|\boldsymbol{\beta}_{\mathcal{S}^c}\|_1 \leq 3\|\boldsymbol{\beta}_{\mathcal{S}}\|_1$,

$$\|\boldsymbol{F}\boldsymbol{\beta}\|^2 \geq n\kappa\|\boldsymbol{\beta}\|^2$$

In our set-up, the *normalized column* condition is automatically met, as $\sup_{x \in \mathcal{X}} |f_j(x)| \leq 1$ for all $j = 1, \ldots, d$ and $x \in \mathcal{X}$. It suffices to check the *restricted eigenvalue* condition holds with high probability. To this end, we adopt the argument (3.3) in Bickel, Ritov and Tsybakov (2009). Recall that $\epsilon_n$ is the element-wise absolute maximum of $\widehat{\boldsymbol{\Sigma}} - \boldsymbol{\Sigma}$. For any $\boldsymbol{\beta}$ such that $\|\boldsymbol{\beta}_{\mathcal{S}^c}\|_1 \leq 3\|\boldsymbol{\beta}_{\mathcal{S}}\|_1$,

$$\frac{\boldsymbol{\beta}'\boldsymbol{F}'\boldsymbol{F}\boldsymbol{\beta}}{n\|\boldsymbol{\beta}\|^2} = \frac{\boldsymbol{\beta}'\boldsymbol{\Sigma}\boldsymbol{\beta}}{\|\boldsymbol{\beta}\|^2} + \frac{\boldsymbol{\beta}'(\widehat{\boldsymbol{\Sigma}} - \boldsymbol{\Sigma})\boldsymbol{\beta}}{\|\boldsymbol{\beta}\|^2} \geq \frac{\boldsymbol{\beta}'\boldsymbol{\Sigma}\boldsymbol{\beta}}{\|\boldsymbol{\beta}\|^2} - \frac{\epsilon_n\|\boldsymbol{\beta}\|_1^2}{\|\boldsymbol{\beta}\|^2}$$
$$\geq \frac{\boldsymbol{\beta}'\boldsymbol{\Sigma}\boldsymbol{\beta}}{\|\boldsymbol{\beta}\|^2} - \epsilon_n\left(\frac{4\|\boldsymbol{\beta}_{\mathcal{S}}\|_1}{\|\boldsymbol{\beta}_{\mathcal{S}}\|}\right)^2 \geq \frac{\boldsymbol{\beta}'\boldsymbol{\Sigma}\boldsymbol{\beta}}{\|\boldsymbol{\beta}\|^2} - 16s\epsilon_n.$$

It follows from Lemma 6.1 that, with probability at least $1 - 2d^{-\delta}$,

$$\epsilon_n \leq \sqrt{2(2+\delta)\alpha(\lambda_{\mathrm{r}} \vee 0)\log d/n}.$$

Assume $\boldsymbol{\Sigma} \succ 0$ and note that $\|\|\boldsymbol{\Sigma}^{-1}\|\|^{-1}$ is the smallest eigenvalue of $\boldsymbol{\Sigma}$. Then the *restricted eigenvalue* condition holds if

$$16\|\|\boldsymbol{\Sigma}^{-1}\|\|s\sqrt{2(2+\delta)\alpha(\lambda_{\mathrm{r}} \vee 0)\log d/n} < 1,$$

and correspondingly

$$\kappa = \|\|\boldsymbol{\Sigma}^{-1}\|\|^{-1} - 16s\sqrt{2(2+\delta)\alpha(\lambda_{\mathrm{r}} \vee 0)\log d/n}.$$

In contrast, the case of i.i.d. data samples requires

$$16\|\|\boldsymbol{\Sigma}^{-1}\|\|s\sqrt{2(2+\delta)\log d/n} < 1,$$

and correspondingly

$$\kappa = \|\|\boldsymbol{\Sigma}^{-1}\|\|^{-1} - 16s\sqrt{2(2+\delta)\log d/n}.$$

The additional multiplicative coefficient $\alpha(\lambda_{\mathrm{r}} \vee 0)$ suggests that the lasso regression on the Markov-dependent data samples requires more data samples than that on i.i.d. data samples and has smaller $\kappa$, which leads to larger error $\|\boldsymbol{\beta}_{w_n} - \boldsymbol{\beta}_\star\|$ and $\|\boldsymbol{\beta}_{w_n} - \boldsymbol{\beta}_\star\|_1$.

6.3. *Sparse covariance.* Let us consider estimating high dimensional covariance matrix with Markov-dependent samples. Suppose that $X_i$'s are samples from a stationary Markov chain $\{X_i\}_{i\geq 1}$ with invariant measure $\pi$ and right spectral gap $1 - \lambda_{\rm r}$; and, $\boldsymbol{f} = (f_1, \ldots, f_d)'$ is a known $d$-dimensional bounded feature mapping. Without loss of generality, we assume $\pi(f_j) = 0$ and $\sup_{x \in \mathcal{X}} |f_j(x)| \leq 1$ for all $j = 1, \ldots, d$. We are interested in estimating the covariance matrix

$$\boldsymbol{\Sigma} = \int \boldsymbol{f}(x)\boldsymbol{f}(x)'\pi(dx).$$

in the high-dimensional regime with $d \gg n$. Define the uniformity class of sparse covariance matrices invariant under permutations by

$$\mathcal{M}(s,m) = \left\{\boldsymbol{M} \succeq 0: \ \boldsymbol{M}_{jj} \leq m, \sum_{k=1}^{d} \mathbb{I}(\boldsymbol{M}_{jk} \neq 0) \leq s, \forall j = 1, \ldots, d\right\}.$$

For a matrix $\boldsymbol{M}$, we define the element-wise thresholding operator by

$$T_t(\boldsymbol{M}) = \left[\boldsymbol{M}_{jk}\mathbb{I}(|\boldsymbol{M}_{jk}| > t)\right]_{1 \leq j,k \leq d}.$$

See Bickel and Levina (2008); Lam and Fan (2009); Cai and Liu (2011); Fan, Liao and Mincheva (2013) for a perspective on this problem. We are ready to present the main theorem of this section. The proof is deferred to Subsection B.1 in the appendix.

THEOREM 6.1. *Suppose $\boldsymbol{\Sigma} \in \mathcal{M}(s,m)$. Let $\widehat{\boldsymbol{\Sigma}} = \sum_{i=1}^{n} \boldsymbol{f}(X_i)\boldsymbol{f}(X_i)'/n$. Take $t$ such that*

$$2\sqrt{2(2+\delta)\alpha(\lambda_r \vee 0)\log d/n} \leq t \lesssim \sqrt{\alpha(\lambda_r \vee 0)\log d/n}.$$

*Then with probability at least $1 - 2d^{-\delta}$*

$$\|\!|T_t(\widehat{\boldsymbol{\Sigma}}) - \boldsymbol{\Sigma}|\!\| \leq \|\!|T_t(\widehat{\boldsymbol{\Sigma}}) - \boldsymbol{\Sigma}|\!\|_1 \leq s\left(2t + 3\sqrt{2(2+\delta)\alpha(\lambda_r \vee 0)\log d/n}\right)$$
$$\lesssim s\sqrt{\alpha(\lambda_r \vee 0)\log d/n}.$$

6.4. *Nonasymptotic error bound for MCMC estimates.* Markov chain Monte Carlo (MCMC) methods have been widely used to compute complicated integrals in Bayesian statistics, machine learning, computational biology, computational physics and computational linguistics. Suppose the task is to compute the integral $\pi(f) = \int f(x)\pi(x)dx$ of a function $f$ with respect to a probability density function $\pi$, but $\pi$ cannot be directly sampled. An MCMC algorithm generates a Markov chain $\{X_i\}_{i\geq 1}$, which converges to the invariant distribution $\pi$, and estimates $\pi(f)$ by averaging the realized values of function $f$ on $n$ MCMC samples $\{X_i\}_{i=n_0+1}^{n_0+n}$ after a burn-in period of length $n_0$.

We have the following bound for the MCMC estimate when $f$ is bounded. The proof is deferred to Subsection B.2 in the appendix.

THEOREM 6.2. *Let $\{X_i\}_{i\geq 1}$ be a Markov chain with invariant measure $\pi$, absolute spectral gap $1 - \lambda$ and right spectral gap $1 - \lambda_r$. Suppose the initial measure $\nu$ is absolutely continuous with respect to the invariant measure $\pi$ and its density, denoted by $d\nu/d\pi$, has a finite p-moment for some $p \in (1, \infty]$. Let $q = p/(p-1) \in [1, \infty)$ and*

$$C_p(\nu, n_0) = \begin{cases} 1 + 2^{2/p}\lambda^{2n_0/q}\left\|\frac{d\nu}{d\pi} - 1\right\|_{\pi,p} & \text{if } p \in (1, 2), \\ 1 + \lambda^{n_0}\left\|\frac{d\nu}{d\pi} - 1\right\|_{\pi} & \text{if } p = 2, \\ 1 + 2^{2/q}\lambda^{2n_0/p}\left\|\frac{d\nu}{d\pi} - 1\right\|_{\pi,p} & \text{if } p \in (2, \infty), \\ \|\frac{d\nu}{d\pi}\|_{\pi,\infty} = \operatorname{ess\,sup} \frac{d\nu}{d\pi} & \text{if } p = \infty. \end{cases}$$

*Then, for any bounded function $f : \mathcal{X} \to [a, b]$ and any $t \in \mathbb{R}$,*

$$\mathbb{E}_\nu\left[e^{t\sum_{i=n_0+1}^{n_0+n} f(X_i) - n\pi(f)}\right] \leq C_p(\nu, n_0)\exp\left(\frac{t^2}{2} \cdot q\alpha(\lambda_r \vee 0) \cdot \frac{n(b-a)^2}{4}\right).$$

*It follows that for any $\epsilon > 0$*

$$\mathbb{P}_\nu\left(\left|\frac{1}{n}\sum_{i=n_0+1}^{n_0+n} f(X_i) - \pi(f)\right| > \epsilon\right) \leq 2C_p(\nu, n_0)\exp\left(-\frac{q^{-1}\alpha(\lambda_r \vee 0)^{-1}}{2 \times (b-a)^2/4} \cdot n\epsilon^2\right).$$

There are mainly two factors influencing the approximation error of MCMC methods. The first factor is how fast the Markov chain converges from the initial measure $\nu$ to the invariant measure $\pi$. This is characterized as $C_p(\nu, n_0)$ in the proceeding theorem. Here $C_p(\nu, n_0)$ upper bounds $\|d(\nu P^{n_0})/d\pi\|_{\pi,p}$ and converges to 1 as $n_0 \to \infty$. The second factor is how the average of $f(X_i)$ fluctuates after the Markov chain reaches the stationarity. This has been already studied by Theorem 2.2.

6.5. *Respondence-driven sampling.* Data about disease prevalence and risk behaviors within some specific subpopulations are crucial in epidemiological studies, but these subpopulations are relatively small and often desire to remain anonymous, rendering standard sampling methods unsuitable. For example, HIV infections are concentrated in three hidden subpopulations: men who have sex with men, injection drug users, and sex workers and their sexual partners (WHO and UNAIDS, 2009). Respondence-driven sampling (RDS), a method initially developed by Heckathorn (1997, 2002) as part of NIH/NIDA-funded HIV prevention project, has now been widely used. RDS collects data through a chain-referral mechanism, in which current participants recruit their contacts to be new participants.

Goel and Salganik (2009) modeled this chain-referral mechanism as a random walk among the network of the population under study (e.g. injection drug users in New York City), and proposed an MCMC importance sampling estimator for the disease prevalence. Formally, let $\mathcal{X}$ denote the population under study, and let $\mathcal{E}$ denote the edge set among them. Let $f(x) = 1$ if a member of the population $x$ is infected and 0 otherwise, let $d(x)$ be the degree of a member $x$ in the network.

Let $\{X_i\}_{i\geq 1}$ be a random walk on the network $(\mathcal{X}, \mathcal{E})$ with uniform edge weights (if assuming participants recruit their contacts uniformly at random). The prevalence of a disease

$$\frac{1}{|\mathcal{X}|} \sum_{x \in \mathcal{X}} f(x)$$

is consistently estimated by

$$\sum_{i=1}^{n} \frac{f(X_i)}{d(X_i)} \bigg/ \sum_{i=1}^{n} \frac{1}{d(X_i)}.$$

Indeed, the random walk has the invariant distribution $\pi(x) = d(x)/2|\mathcal{E}|$, and thus $n^{-1} \sum_{i=1}^{n} f(X_i)/d(X_i) \to \sum_{x \in \mathcal{X}} f(x)/2|\mathcal{E}|$ and $n^{-1} \sum_{i=1}^{n} 1/d(X_i) \to |\mathcal{X}|/2|\mathcal{E}|$.

A random walk on a finite graph is a reversible Markov chain with a self-adjoint $P$. If $P$ admits a right spectral gap $1 - \lambda_{\mathrm{r}}$, then Theorem 6.2 provides a non-asymptotic error bound for the estimator by combining bounds

$$\mathbb{P}_\pi \left( \left| \frac{1}{n} \sum_{i=1}^{n} \frac{f(X_i)}{d(X_i)} - \frac{\sum_{x \in \mathcal{X}} f(x)}{2|\mathcal{E}|} \right| > \epsilon \right) \leq 2 \exp\left( -\frac{2n\epsilon^2}{\alpha(\lambda_{\mathrm{r}} \vee 0)} \right)$$

and

$$\mathbb{P}_\pi \left( \left| \frac{1}{n} \sum_{i=1}^{n} \frac{1}{d(X_i)} - \frac{|\mathcal{X}|}{2|\mathcal{E}|} \right| > \epsilon \right) \leq 2 \exp\left( -\frac{2n\epsilon^2}{\alpha(\lambda_{\mathrm{r}} \vee 0)} \right).$$

6.6. *Multi-armed bandit with Markovian rewards.* The Multi-armed Bandit (MAB) problem has received much attention in decision theory and machine learning. In this problem, there are a number, say $K$ of alternative arms, each with a stochastic reward with initially unknown expectation. The goal is to find the optimal strategy that maximizes the sum of rewards. Let $Z_j(t)$ be the reward from arm $j$ played at round $t$. Let $j_\star(t) = \arg\max_{j=1}^{K} \mathbb{E} Z_j(t)$ be the index of the arm with highest expected reward at round $t$. Let $j(t) \in \{1, \ldots, K\}$ be the arm that is chosen at round $t$. A large body of literatures focuses on minimization of the pseudo-regret on the first $T$ rounds

$$\mathcal{R} = \mathbb{E}\left[ \sum_{t=1}^{T} Z_{j_\star(t)}(t) - \sum_{t=1}^{T} Z_{j(t)}(t) \right].$$

Machine learners have recognized three fundamental formalizations of MAB problems depending on the nature of the reward process: i.i.d., adversarial, and Markovian, but the Markovian formalization has been much less studied than the other two (Bubeck and Cesa-Bianchi, 2012). It is primarily because many bandit algorithms essentially rely on the concentration of average reward draws around its expectation to identify the optimal arm, but powerful concentration inequalities for Markov-dependent random variables are sparse. Let us exemplify the utilities of our theorems by deriving bounds for the pseudo-regret of the celebrated Upper Confidence Bound (UCB) algorithm in the following Markovian MAB problem.

Suppose, in an MAB with $K$ arms, each arm $j$ has an underlying stationary Markov chain $\{X_{ji}\}_{i\geq 1}$ with invariant measure $\pi_j$ and right spectral gap $1-\lambda_{\rm r}$, and a reward function $f_j: x \mapsto [0,1]$. Whenever arm $j$ is played, it returns a reward of the current state and let its Markov chain transition one step. Let $N_j(s) = \sum_{t=1}^{s} \mathbb{I}(j(t)=j)$ be the number of times arm $j$ is played on the first $s$ rounds, and $\Delta_j = \pi_j(f_j) - \pi_{j_\star}(f_{j_\star})$ be the gap between expected rewards of a suboptimal arm $j$ and the optimal arm $j_\star$. Then the pseudo-regret in this set-up is given by

$$\mathcal{R} = \sum_{j=1}^{K} \Delta_j \mathbb{E} N_j(T).$$

Theorem 6.3 bounds the regret for the celebrated UCB algorithm for the MAB problem with Markovian rewards. The proof is deferred to Subsection B.3 in the appendix.

THEOREM 6.3 (*c-UCB algorithm for Markovian MAB*). *Consider the c-UCB algorithm with input parameter $c$ for the above Markovian MAB as follows. Let $\widehat{f}_{j,n} = \frac{1}{n}\sum_{i=1}^{n} f_j(X_{ji})$ be the sample mean of the first $n$ rewards from arm $j$. At each round $t$, select*

$$j(t) \in \arg\max_{j=1}^{K} \widehat{f}_{j,N_j(t-1)} + \sqrt{\frac{c\log t}{2N_j(t-1)}}.$$

*If $c > 2\alpha(\lambda_r \vee 0)$ then this c-UCB algorithm has pseudo-regret*

$$\mathcal{R} \leq \sum_{j:\ \Delta_j > 0} \left( \frac{2c}{\Delta_j} \log T + \frac{c\Delta_j}{c - 2\alpha(\lambda_r \vee 0)} \right).$$

If $\lambda_{\rm r} = 0$ then the proceeding theorem recovers the classical regret bound for $c$-UCB algorithm.

## APPENDIX A: PROOF OF MAIN RESULTS

### A.1. Proof of Theorem 2.3.

PROOF. Changing measure from $\nu$ to $\pi$ and applying Hölder's inequality yields

$$\mathbb{E}_\nu \left[ e^{t\sum_{i=1}^{n} f_i(X_i)} \right] = \mathbb{E}_\pi \left[ \frac{d\nu}{d\pi}(X_1) \cdot e^{t\sum_{i=1}^{n} f_i(X_i)} \right]$$

$$\leq \left\{ \mathbb{E}_\pi \left[ \left|\frac{d\nu}{d\pi}(X_1)\right|^p \right] \right\}^{1/p} \cdot \left\{ \mathbb{E}_\pi \left[ e^{qt\sum_{i=1}^{n} f_i(X_i)} \right] \right\}^{1/q}$$

$$= \left\| \frac{d\nu}{d\pi} \right\|_{\pi,p} \cdot \left\{ \mathbb{E}_\pi \left[ e^{qt\sum_{i=1}^{n} f_i(X_i)} \right] \right\}^{1/q}.$$

Substituting $t$ with $qt$ in (2.2) and plugging it into the above equation yields the first desired inequality. The second desired inequality follows by the Chernoff approach. □

### A.2. Proof of Theorem 2.4.

PROOF. We construct the Markov chain $\{X_i\}_{i\geq 0}$ in the following way. Let $B_1 = 1$ and $\{B_i\}_{i\geq 2}$ be i.i.d. Bernoulli$(1-\lambda)$ random variables, and let $\{W_i\}_{i\geq 1}$ be i.i.d. $\mathcal{N}(0,1)$ random variables. It is easy to verify the construction

$$X_i = (1-B_i)X_{i-1} + B_i W_i, \ \forall i \geq 1.$$

By induction,

$$X_i = \sum_{j=1}^{i} \left( \prod_{k=j+1}^{i} (1-B_k) \right) B_j W_j.$$

Let

$$N_i = \sum_{j=i}^{n} \left( \prod_{k=j+1}^{i} (1-B_k) \right) B_j,$$

then

$$\sum_{i=1}^{n} X_i = \sum_{i=1}^{n} N_i W_i, \quad N_i \geq 0, \quad \sum_{i=1}^{n} N_i = n.$$

Further,

$$\begin{aligned}
\mathbb{E}_\pi \left[ e^{t \sum_{i=1}^n X_i} \right] &= \mathbb{E}_\pi \left[ \mathbb{E}_\pi \left( e^{t \sum_{i=1}^n N_i W_i} \,\Big|\, N_1, \ldots, N_n \right) \right] \\
&= \mathbb{E}_\pi \left[ e^{t^2 \sum_{i=1}^n N_i^2 / 2} \right] \\
&\geq \mathbb{P}_\pi (N_1 = n) \mathbb{E}_\pi \left[ e^{t^2 \sum_{i=1}^n N_i^2 / 2} \,\Big|\, N_1 = n \right] \\
&= \lambda^n e^{t^2 n^2 / 2},
\end{aligned}$$

which could not be bounded by $e^{O(n)t^2/2}$ uniformly for $n \geq 1$ and $t \in \mathbb{R}$. □

### A.3. Results for Markov chains with non-zero spectral gap.

Some literatures such as (Kontoyiannis and Meyn, 2012) define the spectral gap as $1-\lambda_\infty(P)$, where $\lambda_\infty(P)$ is the spectral radius of $P$ acting on $\mathcal{L}_2^0(\pi)$, and coincides with

$$\lambda_\infty(P) = \lim_{k \to \infty} \|\!|P^k - \Pi|\!\|_\pi^{1/k}.$$

Noting that $P\Pi = \Pi P = \Pi$, $\Pi^2 = \Pi$, we can derive that

$$\lambda_\infty(P) = \lim_{k \to \infty} \|\!|(P-\Pi)^k|\!\|_\pi^{1/k} \leq \lim_{k \to \infty} \left[ \|\!|P - \Pi|\!\|_\pi^k \right]^{1/k} = \lambda(P).$$

Thus $1 - \lambda_\infty(P) > 0$ is weaker than $1 - \lambda(P) > 0$. Under this weaker condition, we have Theorem A.1. The idea is to break down $\{X_i\}_{i=1}^n$ into $k$ groups $\{X_i : i \in \mathcal{I}_j\}$ with index sets $\mathcal{I}_j = \{i : 1 \leq i \leq n, (i-1) \equiv j-1 \mod k\}$ for $j = 1, \ldots, k$. Each group forms a Markov chain driven by the $k$-step transition kernel $P^k$ such that

$$\|\!|P^k - \Pi|\!\|_\pi^{1/k} < 1,$$

for a certain $k$, and thus Theorem 2.1 can apply.

THEOREM A.1. *If the Markov operator $P$ of the Markov chain $\{X_i\}_{i\geq 1}$ admits $\lambda_\infty(P) < 1$. Then there exists some $k$ such that $\|\!|P^k - \Pi\|\!|_\pi^{1/k} < 1$. For any bounded functions $f_i : \mathcal{X} \to [a_i, b_i]$ and any $t \in \mathbb{R}$, we have*

$$\mathbb{E}_\pi\left[e^{t\sum_{i=1}^n (f_i(X_i) - \pi(f_i))}\right] \leq \exp\left(\frac{t^2}{2} \cdot k\alpha(\lambda_k^k) \cdot k \max_{j=1}^k \sum_{i \in \mathbb{I}_j} \frac{(b_i - a_i)^2}{4}\right).$$

*The corresponding concentration inequality follows by the Chernoff approach.*

PROOF. Without loss of generality, we assume $\pi(f_i) = 0$. Write

$$\begin{aligned}
\mathbb{E}_\pi\left[e^{t\sum_{i=1}^n f_i(X_i)}\right] &= \mathbb{E}_\pi\left[e^{kt \cdot \frac{1}{k}\sum_{j=1}^k \sum_{i \in \mathcal{I}_j} f_i(X_i)}\right] \\
&\leq \frac{1}{k}\sum_{j=1}^k \mathbb{E}_\pi\left[e^{kt \sum_{i \in \mathcal{I}_j} f_i(X_i)}\right] \quad \text{[Jensen's inequality]} \\
&\leq \frac{1}{k}\sum_{j=1}^k \exp\left(\frac{k^2 t^2}{2} \cdot \alpha(\lambda_k^k) \cdot \sum_{i \in \mathcal{I}_j} \frac{(b_i - a_i)^2}{4}\right) \quad \text{[Theorem 2.1]} \\
&\leq \exp\left(\frac{k^2 t^2}{2} \cdot \alpha(\lambda_k^k) \cdot \max_{j=1}^k \sum_{i \in \mathcal{I}_j} \frac{(b_i - a_i)^2}{4}\right) \\
&= \exp\left(\frac{t^2}{2} \cdot k\alpha(\lambda_k^k) \cdot k \max_{j=1}^k \sum_{i \in \mathcal{I}_j} \frac{(b_i - a_i)^2}{4}\right).
\end{aligned}$$

□

The term $k \max_{j=1}^k \sum_{i \in \mathbb{I}_j} (b_i - a_i)^2/4$ has the same order with $\sum_{i=1}^n (b_i - a_i)/4$, but the term $k\alpha(\lambda_k^k)$ is larger than $\alpha(\lambda_k)$ (Geyer, 1992, Theorem 3.3).

**A.4. Proof of Lemma 4.1.**

PROOF. (i) It is trivial for $g \equiv 0$. For any non-zero $g$,

$$\begin{aligned}
\|g\|_\pi^2 \|\!|M_g \widehat{P}_c M_g\|\!|_\pi &\geq \langle g, M_g \widehat{P}_c M_g g\rangle_\pi \\
&\geq \langle g^2, \widehat{P}_c g^2\rangle_\pi \\
&= c\|g^2\|_\pi^2 + \pi(g^2)^2 \\
&\geq \pi(g^2)^2 = \|g\|_\pi^4.
\end{aligned}$$

Dividing both sides by $\|g\|_\pi^2$ completes the proof.

(ii) Using the facts that $P\Pi = \Pi P = \Pi$ and that $I - \Pi$ is self-adjoint,

$$\begin{aligned}|\langle Ph_1, h_2\rangle_\pi| &= |\langle (I-\Pi)(P-\Pi)(I-\Pi)h_1, h_2\rangle_\pi + \langle \Pi h_1, h_2\rangle_\pi| \\ &= |\langle (P-\Pi)(I-\Pi)h_1, (I-\Pi)h_2\rangle_\pi + \langle \Pi h_1, h_2\rangle_\pi| \\ &\leq |\langle (P-\Pi)(I-\Pi)h_1, (I-\Pi)h_2\rangle_\pi| + |\langle \Pi h_1, h_2\rangle_\pi| \\ &\leq \lambda \|(I-\Pi)h_1\|_\pi \|(I-\Pi)h_2\|_\pi + |\pi(h_1)\pi(h_2)| \\ &\leq \sqrt{\lambda\|(I-\Pi)h_1\|_\pi^2 + \pi(h_1)^2} \cdot \sqrt{\lambda\|(I-\Pi)h_2\|_\pi^2 + \pi(h_2)^2} \\ &= \langle \widehat{P}h_1, h_1\rangle_\pi^{1/2} \cdot \langle \widehat{P}h_2, h_2\rangle_\pi^{1/2}.\end{aligned}$$

(iii) Using part (ii) and the self-adjointness of $S_1$, $S_2$, $S_1\widehat{P}S_1$ and $S_2\widehat{P}S_2$,

$$\begin{aligned}\|\|S_1 P S_2\|\|_\pi &= \sup_{h_1, h_2: \|h_1\|_\pi = \|h_2\|_\pi = 1} |\langle S_1 P S_2 h_2, h_1\rangle_\pi| \\ &= \sup_{h_1, h_2: \|h_1\|_\pi = \|h_2\|_\pi = 1} |\langle P S_2 h_2, S_1 h_1\rangle_\pi| \\ &\leq \sup_{h_1, h_2: \|h_1\|_\pi = \|h_2\|_\pi = 1} \langle \widehat{P} S_1 h_1, S_1 h_1\rangle_\pi^{1/2} \cdot \langle \widehat{P} S_2 h_2, S_2 h_2\rangle_\pi^{1/2} \\ &= \sup_{h_1, h_2: \|h_1\|_\pi = \|h_2\|_\pi = 1} \langle S_1 \widehat{P} S_1 h_1, h_1\rangle_\pi^{1/2} \cdot \langle S_2 \widehat{P} S_2 h_2, h_2\rangle_\pi^{1/2} \\ &= \sup_{h_1: \|h_1\|_\pi = 1} \langle S_1 \widehat{P} S_1 h_1, h_1\rangle_\pi^{1/2} \cdot \sup_{h_2: \|h_2\|_\pi = 1} \langle S_2 \widehat{P} S_2 h_2, h_2\rangle_\pi^{1/2} \\ &= \|\|S_1 \widehat{P} S_1\|\|_\pi^{1/2} \|\|S_2 \widehat{P} S_2\|\|_\pi^{1/2}.\end{aligned}$$

$\square$

## APPENDIX B: PROOF FOR EXAMPLES

### B.1. Proof of Theorem 6.1 for Sparse Covariance.

PROOF. Let $\epsilon_n = \max_{i,j} |\widehat{\boldsymbol{\Sigma}}_{ij} - \boldsymbol{\Sigma}_{ij}|$. Write

$$\begin{aligned}\|\|T_t(\widehat{\boldsymbol{\Sigma}}) - \boldsymbol{\Sigma}\|\|_1 &= \max_i \sum_j |T_t(\widehat{\boldsymbol{\Sigma}}_{ij}) - \boldsymbol{\Sigma}_{ij}| \\ &\leq \max_i \sum_j |\boldsymbol{\Sigma}_{ij}| \mathbb{I}\left(|\widehat{\boldsymbol{\Sigma}}_{ij}| \leq t, |\boldsymbol{\Sigma}_{ij}| \leq t\right) \\ &\quad + \max_i \sum_j |\boldsymbol{\Sigma}_{ij}| \mathbb{I}\left(|\widehat{\boldsymbol{\Sigma}}_{ij}| \leq t, |\boldsymbol{\Sigma}_{ij}| > t\right) \\ &\quad + \max_i \sum_j |\widehat{\boldsymbol{\Sigma}}_{ij} - \boldsymbol{\Sigma}_{ij}| \mathbb{I}\left(|\widehat{\boldsymbol{\Sigma}}_{ij}| > t, |\boldsymbol{\Sigma}_{ij}| > t\right) \\ &\quad + \max_i \sum_j |\widehat{\boldsymbol{\Sigma}}_{ij} - \boldsymbol{\Sigma}_{ij}| \mathbb{I}\left(|\widehat{\boldsymbol{\Sigma}}_{ij}| > t, |\boldsymbol{\Sigma}_{ij}| \leq t\right).\end{aligned}$$

The first term is bounded by $\max_i \sum_{j=1}^d t\mathbb{I}(\boldsymbol{\Sigma}_{ij} \neq 0) \leq ts$; the second term is bounded by $\max_i \sum_{j=1}^d (t+\epsilon_n)\mathbb{I}(\boldsymbol{\Sigma}_{ij} \neq 0) \leq (t+\epsilon_n)s$; the third term is bounded by $\max_i \sum_{j=1}^d \epsilon_n \mathbb{I}(\boldsymbol{\Sigma}_{ij} \neq 0) \leq \epsilon_n s$; and the fourth term is bounded by

$$\max_i \sum_{j=1}^d \epsilon_n \mathbb{I}\left(|\widehat{\boldsymbol{\Sigma}}_{ij}| > t, |\boldsymbol{\Sigma}_{ij}| \leq t\right)$$

$$\leq \epsilon_n \max_i \sum_j \mathbb{I}\left(|\widehat{\boldsymbol{\Sigma}}_{ij}| > t, t/2 < |\boldsymbol{\Sigma}_{ij}| \leq t\right) + \epsilon_n \max_i \sum_j \mathbb{I}\left(|\widehat{\boldsymbol{\Sigma}}_{ij}| > t, |\boldsymbol{\Sigma}_{ij}| \leq t/2\right)$$

$$\leq \epsilon_n s + \epsilon_n \max_i \sum_j \mathbb{I}\left(|\widehat{\boldsymbol{\Sigma}}_{ij} - \boldsymbol{\Sigma}_{ij}| > t/2\right)$$

$$\leq \epsilon_n s + \epsilon_n \max_i \sum_j \mathbb{I}\left(\epsilon_n > t/2\right).$$

Collecting these pieces together yields that if $\epsilon_n \leq t/2$ then

$$|||T_t(\widehat{\boldsymbol{\Sigma}}) - \boldsymbol{\Sigma}|||_1 \leq s(2t + 3\epsilon_n).$$

From Lemma 6.1 it follows that, with probability at least $1 - 2d^{-\delta}$,

$$\epsilon_n \leq \sqrt{2(2+\delta)\alpha(\lambda_r \vee 0)\log d/n}.$$

Thus, setting $t \geq 2\sqrt{2(2+\delta)\alpha(\lambda_r \vee 0)\log d/n}$ yields that, with probability at least $1 - 2d^{-\delta}$,

$$|||T_t(\widehat{\boldsymbol{\Sigma}}) - \boldsymbol{\Sigma}|||_1 \leq s\left(2t + 3\sqrt{2(2+\delta)\alpha(\lambda_r \vee 0)\log d/n}\right).$$

$\square$

**B.2. Proof of Theorem 6.2 for MCMC Estimation.**

PROOF. Let $\nu P^{n_0}$ denote the $n_0$-step transition of $\nu$. Write

$$\mathbb{E}_\nu \left[e^{t\sum_{i=n_0+1}^{n_0+n} f(X_i) - n\pi(f)}\right]$$

$$= \mathbb{E}_{\nu P^{n_0}} \left[e^{t\sum_{i=1}^n f(X_i) - n\pi(f)}\right] \qquad \text{[Markov property]}$$

$$= \mathbb{E}_\pi \left[\frac{d(\nu P^{n_0})}{d\pi}(X_1) \cdot e^{t\sum_{i=1}^n f(X_i) - n\pi(f)}\right] \qquad \text{[Change measure]}$$

$$\leq \left\|\frac{d(\nu P^{n_0})}{d\pi}\right\|_{\pi,p} \times \left\{\mathbb{E}_\pi \left[e^{qt\sum_{i=1}^n f(X_i) - n\pi(f)}\right]\right\}^{1/q} \qquad \text{[Hölder's inequality]}$$

$$\leq \left\|\frac{d(\nu P^{n_0})}{d\pi}\right\|_{\pi,p} \times \exp\left(\frac{t^2}{2} \cdot q\alpha(\lambda_r \vee 0) \cdot \frac{n(b-a)^2}{4}\right). \qquad \text{[Theorem 2.2]}$$

It remains to show $\|d(\nu P^{n_0})/d\pi\|_{\pi,p} \leq C_p(\nu, n_0)$. For $p \in (1, \infty)$, by (Rudolf, 2012, Lemma 3.16 and Proposition 3.17) and the fact that $\lambda(P) = \lambda(P^*)$,

$$\left\|\frac{d(\nu P^{n_0})}{d\pi}\right\|_{\pi,p} \leq \left\|\frac{d(\nu P^{n_0} - \pi)}{d\pi}\right\|_{\pi,p} + 1$$

$$= \left\|[(P^*)^{n_0} - \Pi]\left(\frac{d\nu}{d\pi}\right)\right\|_{\pi,p} + 1 \leq C_p(\nu, n_0).$$

For $p = \infty$, we have

$$\left\|\frac{d(\nu P^{n_0})}{d\pi}\right\|_{\pi,\infty} = \left\|(P^*)^{n_0}\left(\frac{d\nu}{d\pi}\right)\right\|_{\pi,\infty} = \operatorname{ess\,sup}(P^*)^{n_0}\left|\frac{d\nu}{d\pi}\right| \leq \operatorname{ess\,sup}\left|\frac{d\nu}{d\pi}\right|.$$

$\square$

### B.3. Proof of Theorem 6.3 for Multi-armed Bandit.

PROOF. We first argue by contradiction that at least one of three following events given $j(t) = j$ must be true

$$\mathcal{E}_0(t) = \left\{\widehat{f}_{j_\star, N_{j_\star}(t-1)} + \sqrt{\frac{c \log t}{2N_{j_\star}(t-1)}} \leq \pi_{j_\star}(f_{j_\star})\right\},$$

$$\mathcal{E}_1(t) = \left\{\widehat{f}_{j, N_j(t-1)} - \sqrt{\frac{c \log t}{2N_j(t-1)}} > \pi_j(f_j)\right\},$$

$$\mathcal{E}_2(t) = \left\{N_j(t-1) < \frac{2c \log T}{\Delta_j^2}\right\}.$$

Suppose for the sake of contradiction that all the three events are false. Then

$$\widehat{f}_{j_\star, N_{j_\star}(t-1)} + \sqrt{\frac{c \log t}{2N_{j_\star}(t-1)}} > \pi_{j_\star}(f_{j_\star})$$

$$= \pi_j(f_j) + \Delta_j$$

$$\geq \pi_j(f_j) + \sqrt{\frac{2c \log t}{N_j(t-1)}}$$

$$\geq \widehat{f}_{j, N_j(t-1)} + \sqrt{\frac{c \log t}{2N_j(t-1)}},$$

implying $j$ cannot be selected at round $t$, i.e. $j(t) \neq j$.

For any integer $u \geq 1$

$$u \geq \sum_{t=1}^{T} \mathbb{I}(j(t) = j, N_j(t-1) < u)$$

$$\geq \sum_{t=1}^{u} \mathbb{I}(j(t) = j, N_j(t-1) < u) + \sum_{t=u+1}^{T} \mathbb{I}(j(t) = j, N_j(t-1) < u)$$

$$= \sum_{t=1}^{u} \mathbb{I}(j(t) = j) + \sum_{t=u+1}^{T} \mathbb{I}(j(t) = j, N_j(t-1) < u),$$

which implies

$$N_j(T) = \sum_{t=1}^{T} \mathbb{I}(j(t) = j) = \sum_{t=1}^{u} \mathbb{I}(j(t) = j) + \sum_{t=u+1}^{T} \mathbb{I}(j(t) = j)$$

$$\leq u - \sum_{t=u+1}^{T} \mathbb{I}(j(t) = j, N_j(t-1) < u) + \sum_{t=u+1}^{T} \mathbb{I}(j(t) = j)$$

$$= u + \sum_{t=u+1}^{T} \mathbb{I}(j(t) = j, N_j(t-1) \geq u).$$

In particular, let $u = \lceil 2c \log T / \Delta_j^2 \rceil$ then

$$\mathbb{E} N_j(T) \leq u + \sum_{t=u+1}^{T} \mathbb{P}(j(t) = j, N_j(t-1) \geq u)$$

$$\leq u + \sum_{t=u+1}^{T} \mathbb{P}(j(t) = j, \mathcal{E}_2^c(t))$$

$$\leq u + \sum_{t=u+1}^{T} \mathbb{P}(j(t) = j, \mathcal{E}_0(t) \cup \mathcal{E}_1(t))$$

$$\leq u + \sum_{t=u+1}^{T} \mathbb{P}(\mathcal{E}_0(t)) + \sum_{t=u+1}^{T} \mathbb{P}(\mathcal{E}_1(t)).$$

Proceed to bound $\sum_{t=u+1}^{T} \mathbb{P}(\mathcal{E}_0(t))$ and $\sum_{t=u+1}^{T} \mathbb{P}(\mathcal{E}_1(t))$. By Theorem 2.2,

$$\mathbb{P}(\mathcal{E}_0(t)) \leq \sum_{s=1}^{t} \mathbb{P}\left( \widehat{f}_{j,s} - \sqrt{\frac{c \log t}{2s}} > \pi_j(f_j) \right)$$

$$\leq \sum_{s=1}^{t} \exp\left( -\frac{2s}{\alpha(\lambda_r \vee 0)} \times \frac{c \log t}{2s} \right)$$

$$= t^{-c/\alpha(\lambda_r \vee 0) + 1},$$

thus

$$\sum_{t=u+1}^{T} \mathbb{P}(\mathcal{E}_0(t)) \leq \sum_{t=2}^{\infty} t^{-c/\alpha(\lambda_r \vee 0)+1} \leq \frac{1}{c/\alpha(\lambda_r \vee 0) - 2}.$$

The same argument applies for $\sum_{t=u+1}^{T} \mathbb{P}(\mathcal{E}_1(t))$. It follows that

$$\begin{aligned}
\mathbb{E} N_j(T) &\leq u + \frac{2}{c/\alpha(\lambda_r \vee 0) - 2} \\
&\leq \frac{2c \log T}{\Delta_j^2} + 1 + \frac{2}{c/\alpha(\lambda_r \vee 0) - 2} \\
&= \frac{2c \log T}{\Delta_j^2} + \frac{c/\alpha(\lambda_r \vee 0)}{c/\alpha(\lambda_r \vee 0) - 2},
\end{aligned}$$

which completes the proof. □

Department of ORFE  
Princeton University  
205 Sherred Hall  
Princeton, NJ, USA, 08544  
E-mail: jqfan@princeton.edu  
baij@princeton.edu

Department of Statistical Sciences  
University of Toronto  
100 St. George Street  
Toronto, Ontario, Canada, M5S 3G3  
E-mail: qsun@utstat.toronto.edu